\documentclass{article}%
\usepackage{amsmath}
\usepackage{amsthm}
\usepackage{graphicx}
\usepackage{amssymb}
\usepackage{latexsym}
\usepackage{amsbsy}
\usepackage{amsfonts}%
\newtheorem{theorem}{Theorem}
\newtheorem{lemma}[theorem]{Lemma}
\newtheorem{cor}[theorem]{Corollary}
\newtheorem{prop}[theorem]{Proposition}

\theoremstyle{definition}
\newtheorem{remark}[theorem]{Remark}
\newtheorem{definition}[theorem]{Definition}
\newtheorem{hypothesis}[theorem]{Hypothesis}
\begin{document}

\title{Super-Brownian motion with extra birth\\at one point}
\author{Klaus Fleischmann\thanks{Partly supported by the research program
\textquotedblleft Interacting Systems of High Complexity\textquotedblright\ of
the German Science Foundation.}\\{\small Weierstrass Institute for Applied Analysis and Stochastics} \\{\small Mohrenstr.\ 39} \\{\small D\thinspace--\thinspace10117 Berlin} \\{\small Germany}\smallskip\\{\scriptsize E-mail: fleischm@wias-berlin.de} \\{\scriptsize URL: http://www.wias-berlin.de/\symbol{126}fleischm}\\{\scriptsize Fax: +49-30-20\thinspace44\thinspace975}\bigskip
\and Carl Mueller\thinspace\thanks{Supported in part by an NSA grant and the
graduate program \textquotedblleft Probabilistic Analysis and Stochastic
Processes\textquotedblright\ of the German Science Foundation.}\\{\small Department of Mathematics} \\{\small University of Rochester } \\{\small Rochester, NY 14627} \\{\small USA}\smallskip\\{\scriptsize E-mail: cmlr@troi.cc.rochester.edu} \\{\scriptsize URL: http://www.math.rochester.edu/u/cmlr/}\bigskip}
\date{~}
\maketitle

\footnotetext{\emph{Key words and phrases.} Super-Brownian motion,
measure-valued processes, heat equation, singular potential, one-point
potential.} \footnotetext{AMS 1991 \emph{subject classifications.} Primary
60J80, Secondary 60K35.}\setcounter{page}{0}\thispagestyle{empty}

\vfill

\begin{center}
{\scriptsize version of the WIAS Preprint No.\ 790 of November 26, 2002,\quad
ISSN 0946\thinspace--\thinspace8633}
\end{center}

\vfill\newpage$\qquad$\setcounter{page}{1}

\vfill

\begin{quote}
\textbf{Abstract}\quad A super-Brownian motion in two and three dimensions is
constructed where ``particles'' give birth at a higher rate, if they approach
the origin. Via a log-Laplace approach, the construction is based on Albeverio
et al.\ (1995) who calculated the fundamental solutions of the heat equation
with one-point potential in dimensions less than four.
\end{quote}


\vfill

{\scriptsize
\tableofcontents
}

\vfill\newpage

\section{Introduction\label{intro}}

\subsection{Motivation and background}

Measure-valued branching processes, also called superprocesses, arise
naturally as limits of particle branching Markov processes. There is an
immense literature on this topic; see the expositions \cite{Dawson1993},
\cite{Dynkin1994.crm}, \cite{LeGall1999.LN}, \cite{Etheridge2000}, or
\cite{Perkins2000}, for example. Since these models involve mainly
\textquotedblleft non-interacting particles\textquotedblright, many powerful
tools are available, and many detailed properties of these processes are
known. Building on this success, many probabilists have turned their attention
to more complicated models, many of which are governed by singularities. For
example, in 2 or more dimensions, with probability 1, continuous
super-Brownian motion takes values in the space of measures whose closed
support has Lebesgue measure 0. Nevertheless, in certain situations, pairs of
such processes can kill each other when the corresponding \textquotedblleft
particles\textquotedblright\ meet (see, e.g.,\ \cite{EvansPerkins1994}).

Another example of singular behavior is catalytic branching. Here, the
branching of the \textquotedblleft particles\textquotedblright\ is controlled
by a catalytic measure; the higher the \textquotedblleft
density\textquotedblright\ of this measure, the faster the \textquotedblleft
particles\textquotedblright\ branch or die. This catalytic measure can be
supported on a set of Lebesgue measure 0, as long as it is not a polar set of
Brownian motion. In other words, individual \textquotedblleft
particles\textquotedblright\ must have a positive probability of
\textquotedblleft hitting the measure\textquotedblright. See, e.g.,
\cite{DawsonFleischmann1995.high,Delmas1996,FleischmannKlenke1999.smooth,Klenke2000.abscont}%
.

A further example of singular behavior is mass creation. One could imagine a
\textquotedblleft mass creation measure\textquotedblright, which would give
rise to new \textquotedblleft particles\textquotedblright\ whenever the
\textquotedblleft particles\textquotedblright\ of the superprocess hit the
support of the measure. For the extreme case of a single point source
$\delta_{0}$ in $\mathsf{R},$ see
\cite{EnglaenderFleischmann2000.prod,EnglanderTuraev2000}. In particular, a
continuous super-Brownian motion in $\mathsf{R}$ with a point source makes
sense. In higher dimensions, however, at first sight one would expect that a
super-Brownian motion with single point mass creation degenerates to ordinary
super-Brownian motion, since the Brownian particles do not hit a given point.

Our goal is to disprove this intuition. But first, we mention a deterministic,
``one-particle model'', which already gives a different picture. This model
was developed by mathematical physicists starting in the 1930's; see Albeverio
et.\ al.\ \cite{AlbeverioGesztesyHoeg-KrohnHolden1988} for historical background.

Consider the heat equation in $\mathsf{R}^{d}$ with a one-point potential:
\begin{equation}
\frac{\partial u}{\partial t}\;=\;\Delta u\,+\,\delta_{0}^{(\alpha
)}u\;=:\;\Delta^{\!(\alpha)}u \label{singular}%
\end{equation}
(by $f:=g$ or $g=:f$, we mean that $f$ is defined to be equal to $g$).
Heuristically, $\Delta^{\!(\alpha)}=\Delta+\delta_{0}^{(\alpha)}$ is the limit
as $\varepsilon\downarrow0$ of the operator $\,$%
\begin{equation}
\Delta_{\varepsilon}^{\!(\alpha)}\ :=\ \Delta\,+\,h(d,\alpha,\varepsilon
)\,\varepsilon^{-d}\mathbf{1}_{B_{\varepsilon}(0)\,},
\end{equation}
where $B_{\varepsilon}(y)$ denotes the open ball around $y\in\mathsf{R}^{d}$
with radius $\,\varepsilon,$\thinspace\ and $\,h(d,\alpha,\varepsilon
)$\thinspace\ is some additional rescaling factor, depending on a parameter
$\alpha,$ at least.

For instance, in dimension $d=3$,
\begin{equation}
h(3,\alpha,\varepsilon)\;:=\;\big(k+\tfrac{1}{2}\big)^{\!2}\pi^{2}%
\varepsilon-8\pi^{2}\alpha\varepsilon^{2}-\zeta\varepsilon^{3},\qquad\alpha
\in\mathsf{R},\quad\varepsilon>0, \label{not.h}%
\end{equation}
where $\,k$\thinspace\ is any integer and $\,\zeta$\thinspace\ any real number
(we rely on \cite[(H.74)]{AlbeverioGesztesyHoeg-KrohnHolden1988}). Then, in a
sense, $\Delta_{\varepsilon}^{\!(\alpha)}\rightarrow\Delta^{\!(\alpha)}$ as
$\varepsilon\downarrow0,$ where the limit operator $\Delta^{\!(\alpha)}$ is
\emph{independent}\/ of $k$ and $\zeta$ (so for simplification one could set
$\,k=0=\zeta).$\thinspace\ 

Actually, in the physics literature, main attention is paid to the
Schr\"{o}dinger equation instead of the heat equation (and the positive
definite operator $-\Delta$ is preferred instead of $\Delta)$. But note that
in the time-stationary case, Schr\"{o}dinger and heat equation coincide. So
the operators $\Delta^{\!(\alpha)}$ are relevant in both cases. Physically,
the parameter $\alpha$ is related to the \textquotedblleft scattering
length\textquotedblright\ $\,(4\pi\alpha)^{-1}$ (see, for instance,
\cite[p.13]{AlbeverioGesztesyHoeg-KrohnHolden1988}). In particular,
$\Delta^{\!(\alpha)}\rightarrow\Delta$ as $\alpha\uparrow\infty,$ giving the
\textquotedblleft free\textquotedblright\ case. We understand $\delta
_{0}^{(\alpha)}$ as $\lambda_{\alpha}\delta_{0\,}$. In dimension
$\,d=3$\thinspace\ the coupling constant $\lambda_{\alpha}$ of the point
source $\,\delta_{0}$ has to be of the form $\lambda_{\alpha}=\varepsilon
-\alpha\varepsilon^{2}$ with $\varepsilon$ \textquotedblleft
infinitesimal\textquotedblright\ in a special way.

Even though the number of particles that hit the origin is infinitesimal, one
can imagine that they give raise to a positive mass, provided that the birth
rate is high enough. This explains, why the linear $\varepsilon$--term in
(\ref{not.h}) is not allowed to be too small, in particular, it cannot be
negative. In the latter case, particles will simply die, and nothing else will
happen. But since there are only infinitesimally many particles hitting the
origin, their possible death will pass unnoticed. 

At this point we would like to understand certain questions from a
probabilistic point of view. For instance, in dimensions $\,d=3$\thinspace,
why don't all sufficiently high coefficients of the linear $\varepsilon$--term
occur, and why is the limit operator $\Delta^{\!(\alpha)}$\thinspace
\ independent of the integer $\,k$?\thinspace\ Unfortunately, this is outside
the scope of the present paper.

Strictly speaking, $\,\{\Delta^{\!(\alpha)}:\,\alpha\in\mathsf{R}%
\}\cup\{\Delta\}\,$ is the \vspace{-1pt}family of all self-adjoint extensions
on $\,\mathcal{L}^{2}(\mathsf{\dot{R}}^{d},\mathrm{d}x),$\thinspace
\ $d=2,3,$\thinspace\ of the Laplacian $\Delta$ acting on $\,\mathcal{C}%
_{0}^{(\infty)}(\mathsf{\dot{R}}^{d}),$\thinspace\ where $\,\mathsf{\dot{R}%
}^{d}:=\mathsf{R}^{d}\backslash\{0\}.$\thinspace\ See, for instance,
\cite[Chapters I.1 and I.5]{AlbeverioGesztesyHoeg-KrohnHolden1988}.

The \emph{fundamental solutions}\/ $P^{\alpha}$ to equation
\begin{equation}
\frac{\partial u}{\partial t}\;=\;\Delta^{\!(\alpha)}u\quad\text{on}%
\quad(0,\infty)\times\mathsf{\dot{R}}^{d},\quad d=2,3,\label{singular'}%
\end{equation}
have been computed in \cite{AlbeverioBrzezniakDabrowski1995}. $P^{\alpha}$ is
different from the heat kernel, for each $\,\alpha\in\mathsf{R}$%
\thinspace\ (only for $\,\alpha\uparrow\infty$\thinspace\ one gets back the
heat kernel; for the case $\,d=3,$\thinspace\ see Subsection \ref{SS.P.alpha}
below). $P^{\alpha}$ is a basic object in the present paper.

\subsection{Sketch of result}

Based on the preceding analytical results from
\cite{AlbeverioBrzezniakDabrowski1995}, the \emph{purpose}\/ of the present
paper is to construct a measure-valued super-Brownian motion $X=\{X_{t}%
:\,t\geq0\}$ in $\,\mathsf{\dot{R}}^{d}=\mathsf{R}^{d}\backslash
\{0\},$\thinspace\ $d=2,3,$\thinspace\footnote{$^{)}$ Recall that the
one-dimensional super-Brownian motion with extra birth at $0,$ that is,
related to the log-Laplace equation%
\[
\frac{\partial v}{\partial t}\;=\;\frac{1}{2}\Delta v\,+\delta_{0}%
\,-\,v^{2}\quad\text{on}\quad(0,\infty)\times\mathsf{R},
\]
was introduced in \cite{EnglaenderFleischmann2000.prod}.}\thinspace\ related
to the formal log-Laplace equation
\begin{equation}
\left.
\begin{array}
[c]{c}%
\dfrac{\partial v}{\partial t}\;=\;\Delta^{\!(\alpha)}v-\eta\,v^{1+\beta}%
\quad\text{on \thinspace}(0,\infty)\times\mathsf{\dot{R}}^{d},\vspace{6pt}\\
v(0+,x)\;=\;\varphi(x)\;\geq\;0,\qquad x\in\mathsf{\dot{R}}^{d},
\end{array}
\right\}  \label{cumulant}%
\end{equation}
with constants $\,0<\beta\leq1,$\thinspace\ $\eta\geq0,$\thinspace\ and where
the $\,\varphi$\thinspace\ are appropriate test functions. Of course, $X$ is
related to (\ref{cumulant}) via the log-Laplace transition functional%
\begin{equation}
-\log\mathbf{P}\big\{\mathrm{e}^{-\left\langle X_{t\,},\varphi\right\rangle
}\;\big|\;X_{0}\big\}\;=\;\left\langle X_{0\,},v(t,\,\cdot\,)_{\!_{\!_{\,}}%
}\right\rangle \!,\qquad t>0, \label{logLap}%
\end{equation}
of the Markov process $\,X.$\thinspace\ 

Roughly speaking, we have ``many'' independent Brownian ``particles'' which
everywhere undergo critical branching with index $1+\beta$ and rate $\eta,$
but additionally give birth to new particles if they ``approach'' the origin
$0.$

Here is a rough formulation of our main result; a more precise statement will
be given in Theorem \ref{T.ex} in Subsection \ref{SS.constrX} below.

\begin{theorem}
[\textbf{Existence of }$X$]\label{T.ex.0}If $\,d=2,$\thinspace\ let
$\,0<\beta\leq1,$\thinspace\ and if $\,d=3,\,\ $let $\,0<\beta<1.$%
\thinspace\ Then, for each $\,\alpha\in\mathsf{R},$\thinspace\ there is a
(unique in law) non-degenerate measure-valued (time-homogeneous) Markov
process $\,X=X^{\alpha}$\thinspace\ having log-Laplace transition functional\/
\emph{(\ref{logLap})} with $\,v$\thinspace\ solving\/ \emph{(\ref{cumulant}).}
\end{theorem}

We call $\,X$\thinspace\ a \emph{super-Brownian motion in}\/ $\,\mathsf{R}%
^{d}$\thinspace\ \emph{with extra birth at point}\/ $\,x=0.$\thinspace\ Note
that in the case $\,\eta=0,$\thinspace\ the process degenerates to the
deterministic mass flow related to the kernels $\,P^{\alpha},$\thinspace\ the
fundamental solutions to (\ref{singular'}). At the same time, this mass flow
is identical with the expectation of $\,X,$\thinspace\ for any $\,\eta
.$\thinspace\ In particular, $\,X=X^{\alpha}$\thinspace\ is different from
ordinary super-Brownian motion (corresponding to $\,\alpha=\infty).$

\begin{remark}
[\textbf{Open problem}]The condition $\,\beta<1$\thinspace\ in the
three-dimensional case, which excludes finite variance branching as in
continuous super-Brownian motion, looks a bit strange. We need this condition
for technical reasons, to handle some singularities at the point
$\,x=0$\thinspace\ where extra birth occurs (see Remark \ref{R.beta}
below).\hfill$\Diamond$
\end{remark}

It would, of course, be interesting to reveal that this superprocess
$\,X$\thinspace\ has strange new properties. However, we leave this task for a
future paper and present this construction result separately, since it seems
to be interesting enough.

\subsection{Outline}

In Section \ref{S.basic} we give some estimates involving the basic solutions
and the semigroup related to the linear equation (\ref{singular'}). Then, in
Section \ref{for-both}, we show that the log-Laplace equation (\ref{cumulant})
is well posed. Here we use Picard iteration, but for the non-negativity of
solutions we go back to a linearized equation. For the construction of
$\,X$\thinspace\ in Section \ref{S.Constr.X} we use a Trotter product formula,
alternating between purely continuous state branching (Feller's branching
diffusion if $\,\beta=1)$ and deterministic mass flow with single point mass
creation (related to the kernels $\,P^{\alpha}).$

For background from a mathematical physics point of view concerning the
operators $\Delta^{\!(\alpha)}$ we recommend
\cite{AlbeverioGesztesyHoeg-KrohnHolden1988}, and for basic facts on
superprocesses we refer to one of the systematic treatments
\cite{Dawson1993,Dynkin1994.crm,LeGall1999.LN,Etheridge2000,Perkins2000} which
we have already mentioned.

\section{The heat equation with birth at a single point\label{S.basic}}

After introducing the set $\,\Phi$\thinspace\ of test functions, on which the
heat flow acts continuously (Lemma \ref{L.cont.heat}), we define the kernel
$\,P^{\alpha}$\thinspace\ in Subsection \ref{SS.P.alpha} and show the strong
continuity of the related flow $\,S^{\alpha}$\thinspace\ on $\,\Phi$%
\thinspace\ (Corollary \ref{C.scSalpha}).

\subsection{Preliminaries: test functions in $\Phi\subset\mathcal{H}_{+}%
$\label{SS.Prel}}

The letter $C$ denotes a constant which might change its value from occurrence
to occurrence. $C_{\#}\,\ $and $C_{(\#)}$ refer to specific constants which
are defined around Lemma $\#,$ say, or formula $(\#),$ respectively.

Let $\phi$ denote the \emph{weight and reference function}%
\begin{equation}
\phi(x)\;:=\;|x|^{-(d-1)/2},\qquad x\in\mathsf{\dot{R}}^{d}=\mathsf{R}%
^{d}\backslash\{0\}. \label{ref.fct}%
\end{equation}

\textrm{F}or each fixed constant $\,\varrho\geq1,$\thinspace\ we introduce the
Lebesgue space $\mathcal{H}=\mathcal{H}^{\varrho}=\mathcal{L}^{\varrho
}\!\big
(\mathsf{\dot{R}}^{d},\,\phi(x)\mathrm{d}x\big)$ of equivalence classes
$\,\varphi$\thinspace\ of measurable functions on $\,\mathsf{\dot{R}}^{d}\,$
for which $\Vert\varphi\Vert_{\mathcal{H}}<\infty$, where\thinspace
\footnote{$^{)}$ Here and in similar cases we use this simplified integration
domain $\,\mathsf{R}^{d}$, since in case of integration with respect Lebesgue
measure there is no difference including or not including the point
$\,0\in\mathsf{R}^{d}$\thinspace\ of singularity.}
\begin{equation}
\left\Vert \varphi\right\Vert _{\mathcal{H}}\;:=\;\Big(\int_{\mathsf{R}^{d}%
}\mathrm{d}x\;\phi(x)\,|\varphi|^{\varrho}(x)\Big)^{\!1/\varrho}%
.\label{norm.H}%
\end{equation}
(As usual, we do not distinguish between an equivalence class and its representatives.)

For fixed $\,\varrho\geq1,$\thinspace\ let $\,\Phi=\Phi^{\varrho}$%
\thinspace\ denote the set of all \emph{continuous}\/ functions $\,\varphi
:\mathsf{\dot{R}}^{d}\rightarrow\mathsf{R}$\thinspace\ such that $\,\varphi
\in\mathcal{H}$\thinspace\ and%
\begin{equation}
0\,\leq\,\varphi\,\leq\,C_{\mathrm{(\ref{phi.bound})}}\phi\quad\text{for some
constant}\,\;C_{\mathrm{(\ref{phi.bound})}}=C_{\mathrm{(\ref{phi.bound})}%
}(\varphi). \label{phi.bound}%
\end{equation}
We endow $\,\Phi$\thinspace\ with the topology inherited from $\,\mathcal{H}%
.$\thinspace\ Note that the set $\,\mathcal{C}_{\mathrm{com}}^{+}%
=\mathcal{C}_{\mathrm{com}}^{+}(\mathsf{\dot{R}}^{d})$\thinspace\ of all
non-negative continuous functions on $\,\mathsf{\dot{R}}^{d}$\thinspace\ with
compact support is contained in $\,\Phi.$\thinspace\ Note also that
$\,\varphi\in\Phi$\thinspace\ might have a singularity at $\,x=0$%
\thinspace\ of order $\,|x|^{-\xi}$\thinspace\ with $\,0<\xi<\frac{d+3}{2}%
.$\thinspace\ The functions in $\,\Phi$\thinspace\ will serve as test
functions in log-Laplace representations.

Let $\,\mathcal{M}=\mathcal{M}(\mathsf{\dot{R}}^{d})$\thinspace\ denote the
set of all measures $\,\mu$\thinspace\ defined on $\,\mathsf{\dot{R}}^{d}%
$\thinspace\ such that $\,\left\langle \mu,\varphi\right\rangle :=\int
\!\mu(\mathrm{d}x)\,\varphi(x)<\infty$\thinspace\ for all $\,\varphi\in\Phi
.$\thinspace\ We equip $\,\mathcal{M}$\thinspace\ with the vague topology
(recall that $\,\mathcal{C}_{\mathrm{com}}^{+}\subset\Phi).$\thinspace\ Of
course, each measure $\,\mu\in\mathcal{M}$\thinspace\ can also be considered
as a measure on $\,\mathsf{R}^{d}$\thinspace\ with zero mass at $\,0\in
\mathsf{R}^{d}.$\thinspace\ But in our pairing $\left\langle \mu
,\varphi\right\rangle ,$\thinspace\ $\varphi\in\Phi,$\thinspace\ we cannot
extend to work with measures $\,\mu$\thinspace\ on $\,\mathsf{R}^{d}%
$\thinspace\ allowing positive mass at $\,0$\thinspace\ by the mentioned
possible singularities of the $\,\varphi\in\Phi.$

If $\,\mu$\thinspace\ is a finite measure, we write $\,\Vert\mu\Vert
$\thinspace\ for its total mass. The symbol $\,\ell$\thinspace\ denotes the
Lebesgue measure, $\,A^{\mathrm{c}}$\thinspace\ the complement of
$\,A,$\thinspace\ and $\,a\vee b$\thinspace\ the maximum of $\,a$%
\thinspace\ and $\,b.$

\subsection{Heat flow estimates on $\mathcal{H}$\label{SS.heat.flow}}

In this subsection we fix a dimension $d\geq1.$ Let $P=P(t;x,y)$ refer to the
fundamental solution of the heat equation $\,$%
\begin{equation}
\frac{\partial u}{\partial t}\;=\;\Delta u\quad\text{on}\quad(0,\infty
)\times\mathsf{R}^{d}.
\end{equation}
In other words,
\begin{equation}
P(t;x,y)\;=\;(4\pi t)^{-d/2}\,\mathrm{e}^{-|y-x|^{2}/4t}\,\!,\qquad t>0,\quad
x,y\in\mathsf{R}^{d}.
\end{equation}
Let $S=\{S_{t}:\,t\geq0\}\ $denote the semigroup corresponding to this heat
kernel $P$.

Here is our first estimate [with $\,\phi\,$ the weight and reference function
introduced in (\ref{ref.fct})]:

\begin{lemma}
[A heat flow estimate]\label{L.hf.est}There is a constant
$\,C_{\mathrm{\ref{L.hf.est}}}=C_{\mathrm{\ref{L.hf.est}}}(d)$\thinspace\ such
that%
\begin{equation}
S_{t}\phi\;\leq\;C_{\mathrm{\ref{L.hf.est}}}\,\phi,\qquad t\geq0.
\end{equation}

\end{lemma}

%

\proof
Without loss of generality, let $\,t>0$\thinspace\ and $\,x\neq0.$%
\thinspace\ We have to show that%
\begin{equation}
\frac{1}{\phi(x)}\,S_{t}\phi\,(x)\;=\;\frac{1}{\phi(x)}\int_{\mathsf{R}^{d}%
}\!\mathrm{d}y\;\phi(y)\,\frac{1}{(4\pi t)^{d/2}}\,\,\mathrm{e}^{-|y-x|^{2}%
/4t} \label{start}%
\end{equation}
is bounded in $t>0$ and $x\neq0.$ By the change of variables $\,y\rightarrow
t^{1/2}y$\thinspace\ and with notation $\,xt^{-1/2}=:z\neq0,$\thinspace\ we
get%
\begin{gather}
\frac{1}{\phi(x)}\,S_{t}\phi(x)\;=\;\frac{C}{\phi(z)}\,\int_{\mathsf{R}^{d}%
}\!\mathrm{d}y\;\phi(y)\,\,\mathrm{e}^{-|y-z|^{2}/4}\nonumber\\[2pt]
\;\leq\;\frac{C}{\phi(z)}\,\int_{\mathsf{R}^{d}}\!\mathrm{d}y\;\phi
(y)\,\,\mathrm{e}^{-\left[  |y|-|z|_{\!_{\!_{\,}}}\right]  ^{2}/4}.
\end{gather}
Setting $\,|z|=:s>0$\thinspace\ and substituting $\,|y|=r,$\thinspace\ we
further obtain%
\begin{equation}
\frac{1}{\phi(x)}\,S_{t}\phi(x)\;\leq\;C\,s^{(d-1)/2}\int_{0}^{\infty
}\!\mathrm{d}r\;r^{(d-1)/2}\,\,\,\mathrm{e}^{-|r-s|^{2}/4}.
\end{equation}
First we restrict the integration to $\,r\leq2s.$\thinspace\ For this part we
get the bound%
\begin{equation}
C\,s^{d-1}\int_{0}^{2s}\!\mathrm{d}r\;\,\,\mathrm{e}^{-|r-s|^{2}%
/4}\,\,\;=\;C\int_{0}^{2}\!\mathrm{d}r\;s^{d}\,\mathrm{e}^{-s^{2}|r-1|^{2}/4},
\end{equation}
by the substitution $r\rightarrow sr.$ The latter expression is continuous in
$s>0$ and converges to $0$ as $s\downarrow0$ and as $s\uparrow\infty$ (use
bounded convergence for $r\neq1),$ thus it is bounded in $s>0.$ On the other
hand, for the integration restricted to $r\geq2s$ we get the bound%
\begin{equation}
C\int_{2s}^{\infty}\!\mathrm{d}r\;r^{d-1}\,\,\mathrm{e}^{-|r-s|^{2}%
/4}\,\,\;=\;C\int_{s}^{\infty}\!\mathrm{d}r\;(r+s)^{d-1}\,\mathrm{e}%
^{-r^{2}/4}.
\end{equation}
But $\,(r+s)/r=1+s/r\leq2$\thinspace\ for $\,r\geq r\geq s>0,$\thinspace\ and%
\begin{equation}
\int_{0}^{\infty}\!\mathrm{d}r\;r^{d-1}\,\mathrm{e}^{-r^{2}/4}\;<\;\infty.
\end{equation}
This completes the proof.\hfill\rule{0.5em}{0.5em}\bigskip

We finish this subsection with a simple maximization result.

\begin{lemma}
[\textbf{Maximum in the center}]\label{L.max}Fix a constant $\,\varkappa
>0.$\thinspace\ Then%
\begin{equation}
S_{t}\phi^{\varkappa}\,(x)\;\leq\;S_{t}\phi^{\varkappa}\,(0),\qquad t>0,\quad
x\in\mathsf{R}^{d}. \label{max}%
\end{equation}

\end{lemma}

%

\proof
We will use the fact that in the integral%
\begin{equation}
\int_{\mathsf{R}^{d}}\mathrm{d}y\;P(t;x,y)\,\phi^{\varkappa}(y) \label{int}%
\end{equation}
the mapping $\,y\mapsto\phi^{\varkappa}(y)$\thinspace\ is radially symmetric
and decreasing in $\,|y|.\,\ $The same is true for $\,y\mapsto P(t;x,y),$%
\thinspace\ except a shift by $\,x.$\thinspace\ \medskip

\noindent1$^{\circ}$ (\emph{Simplification}).\quad Let $\,a,b,c,d\geq
0,$\thinspace\ then, by expanding,%
\begin{equation}
(a+b)(c+d)+ac\;\geq\;(a+b)c+a(c+d).
\end{equation}

\noindent2$^{\circ}$ (\emph{Functions with }$n$ \emph{steps}).\quad For
$\,n\geq2,$\thinspace\ let $\,$%
\begin{equation}
f_{i}\;:=\;\sum_{j=1}^{n}a_{i,j}\,\mathsf{1}_{B_{j}}\;\geq\;0,\qquad
i=1,2,\label{step}%
\end{equation}
be two step functions defined on $\,n\geq2$\thinspace\ cubes $\,B_{1\,}%
,\ldots,B_{n}$\thinspace\ in $\,\,\mathsf{R}^{d}$\thinspace\ of equal volume,
say $\,\mathrm{v}.$\thinspace\ For $\,i=1,2,$\thinspace\ let $\,\bar{f}_{i}%
$\thinspace\ be constructed from $\,f_{i}$\thinspace\ by rearranging the
$\,a_{i,j}$\thinspace\ to $\,\bar{a}_{i,1}\geq\cdots\geq\bar{a}_{i,n\,}%
.$\thinspace\ Then%
\begin{equation}
\int_{\mathsf{R}^{d}}\mathrm{d}x\;\bar{f}_{1}(x)\,\bar{f}_{2}(x)\;\geq
\;\int_{\mathsf{R}^{d}}\mathrm{d}x\;f_{1}(x)\,f_{2}(x).\label{ff}%
\end{equation}
In fact,
\begin{equation}
\int_{\mathsf{R}^{d}}\mathrm{d}x\;f_{1}(x)\,f_{2}(x)\;=\;\mathrm{v}\sum
_{j=1}^{n}a_{1,j}\,a_{2,j}\,.\label{integ}%
\end{equation}
Rearranging if necessary, we may assume that $\,f_{1}=\bar{f}_{1\,}%
,$\thinspace\ that is $\,a_{1,j}=\bar{a}_{1,\,j\,},$\thinspace\ $1\leq j\leq
n.$\thinspace\ Exploiting step 1$^{\circ},$ we may switch from $\,f_{2}\,\ $to
$\,\,\bar{f}_{2}$\thinspace\ by a sequence of rearrangements which never
decrease the integral in (\ref{integ}). This then gives the claim
(\ref{ff}).\medskip

\noindent3$^{\circ}$ (\emph{Approximation}).\quad We may assume that the right
hand side of (\ref{max}) is finite. Then the ``integrals'' in (\ref{max})
[recall (\ref{int})] can be approximated by using step functions as in
(\ref{step}). Then (\ref{max}) follows from (\ref{ff}) by passing to the
limit.\hfill\rule{0.5em}{0.5em}

\subsection{Strong continuity of the heat flow on $\,\mathcal{H}$}

Next we will proof the following statement.

\begin{lemma}
[Estimate of $S$ in case of an additional singularity]\label{new-mink}%
\hspace{-2.1pt}Let $\,d$ $\geq1,\,\ 0\leq\beta\leq1,$\thinspace\ and assume
that $\,\varrho$\thinspace\ in\/ \emph{(\ref{norm.H}) }satisfies%
\begin{equation}
\varrho\;>\;\frac{1}{1-\beta(d-1)/2d}\,. \label{ass}%
\end{equation}
Then there is a constant $\,C_{\mathrm{\ref{new-mink}}\,}%
=C_{\mathrm{\ref{new-mink}}\,}(d,\beta,\varrho)$\thinspace\ such that for all
$\,\varphi\in\mathcal{H}=\mathcal{H}^{\varrho},$%
\begin{equation}
\left\Vert S_{t}(\varphi\phi^{\beta})\right\Vert _{\mathcal{H}}^{\varrho
}\;\leq\;C_{\mathrm{\ref{new-mink}}}\,t^{-\beta\varrho(d-1)/4}\,\Vert
\varphi\Vert_{\mathcal{H}}^{\varrho}\,,\qquad t>0. \label{new-mink'}%
\end{equation}

\end{lemma}

%

\proof
Fix $\,d,\beta,\varrho$\thinspace\ as in the lemma. For $t>0\,$ and
$\,x\in\mathsf{R}^{3},$\thinspace\ we introduce the measures%
\begin{equation}
\mu_{t,x}(\mathrm{d}y)\;:=\;t^{\kappa}P(t;x,y)\,\phi^{\lambda}(y)\,\mathrm{d}y
\end{equation}
with%
\begin{equation}
\kappa\;:=\;\frac{\beta\varrho(d-1)}{4(\varrho-1)}\quad\text{and}\quad
\lambda\;:=\;\frac{\beta\varrho}{\varrho-1}\,. \label{constants}%
\end{equation}
By Lemma \ref{L.max},
\begin{equation}
\Vert\mu_{t,x}\Vert\;\leq\;\Vert\mu_{t,0}\Vert\;=\;\int_{\mathsf{R}^{d}%
}\mathrm{d}y\;P(1;0,y)\,\phi^{\lambda}(y)\;=:\;C_{(\mathrm{\ref{C4}})},
\label{C4}%
\end{equation}
where in the last step we used Brownian scaling and the identity
$\,\kappa-\lambda(d-1)/4=0.$\thinspace\ Note that $\,C_{(\mathrm{\ref{C4}}%
)}=C_{(\mathrm{\ref{C4}})}(d,\beta,\varrho)$\thinspace\ is finite by our
assumption (\ref{ass}). Therefore, the measures $\,\mu_{t,x}$\thinspace\ are
finite with total mass at most $\,C_{(\mathrm{\ref{C4}})}$\thinspace
\ independent of $\,t$ and $x.$ Now, for each finite measure $\,\mu$%
\thinspace\ on $\,\mathsf{R}^{d},$ and measurable $\,\varphi,$\thinspace\ by
H\"{o}lder's inequality,%
\begin{equation}
\Big[\int_{\mathsf{R}^{d}}|\varphi|(y)\,\mu(\mathrm{d}y)\Big]^{\varrho}%
\;\leq\;\Vert\mu\Vert^{\varrho-1}\,\int_{\mathsf{R}^{d}}\mu(\mathrm{d}%
y)\;|\varphi|^{\varrho}(y). \label{Hoelder}%
\end{equation}
Applied to the measures $\,\mu_{t,x}$\thinspace\ we get%
\begin{align}
\left\vert S_{t}(\varphi\phi^{\beta})\,(x)\right\vert ^{\varrho}\;  &
=\;t^{-\kappa\varrho}\Big|\int_{\mathsf{R}^{d}}\mu_{t,x}(\mathrm{d}%
y)\,\phi^{\beta-\lambda}(y)\,\varphi(y)\Big|^{\varrho}\nonumber\\
&  \leq\;t^{-\kappa\varrho}\,\Vert\mu_{t,x}\Vert^{\varrho-1}\int
_{\mathsf{R}^{d}}\mu_{t,x}(\mathrm{d}y)\,\phi^{(\beta-\lambda)\varrho
}(y)\,|\varphi|^{\varrho}(y)\nonumber\\[3pt]
&  \leq\;t^{-\kappa(\varrho-1)}\,C_{(\mathrm{\ref{C4}})}^{\varrho-1}%
\,S_{t}|\varphi|^{\varrho}\,(x),
\end{align}
since $\,\lambda+(\beta-\lambda)\varrho=0$\thinspace\ by (\ref{constants}).
But by Lemma \ref{L.hf.est},%
\begin{gather}
\int_{\mathsf{R}^{d}}\mathrm{d}x\;\phi(x)\,S_{t}|\varphi|^{\varrho
}\,(x)\;=\;\int_{\mathsf{R}^{d}}\!\mathrm{d}y\;|\varphi|^{\varrho}%
(y)\,S_{t}\phi\,(y)\;\nonumber\\[2pt]
\leq\;\int_{\mathsf{R}^{d}}\!\mathrm{d}y\;|\varphi|^{\varrho}%
(y)\,C_{\mathrm{\ref{L.hf.est}}}\,\phi(y)\;=\;C_{\mathrm{\ref{L.hf.est}}%
}\,\Vert\varphi\Vert_{\mathcal{H}}^{\varrho}\,.
\end{gather}
Hence,%
\begin{equation}
\Vert S_{t}(\varphi\phi^{\beta})\Vert_{\mathcal{H}}^{\varrho}\;\leq
\;t^{-\kappa(\varrho-1)}\,C_{(\mathrm{\ref{C4}})}^{\varrho-1}%
\,C_{\mathrm{\ref{L.hf.est}}}\,\Vert\varphi\Vert_{\mathcal{H}}^{\varrho}\,,
\end{equation}
and the claim follows since $\,\kappa(\varrho-1)=\beta\varrho(d-1)/4.$%
\hfill\rule{0.5em}{0.5em}\bigskip

Lemma \ref{new-mink} with $\,\beta=0$\thinspace\ yields the following result.

\begin{lemma}
[\textbf{Strong continuity of the heat flow on }$\mathcal{H}$]%
\label{L.cont.heat}The semigroup $\,S$\thinspace\ acting on $\,\mathcal{H=H}%
^{\varrho}$\thinspace\ is strongly continuous.
\end{lemma}

%

\proof
Fix $\,\varphi\in\mathcal{H}.$\thinspace\ By linearity, we may assume that
$\,\varphi\geq0.$\thinspace\ Consider $\,t\in(0,1].$\medskip

\noindent1$^{\circ}$ (\emph{Reducing to bounded functions on compact sets}).
Fix $\,\varepsilon\in(0,1].$\thinspace\ We choose a compact set $\,K\subset
\mathsf{R}^{d}$\thinspace\ so large that $\,\Vert\varphi\,\mathsf{1}%
_{K^{\mathrm{c}}}\Vert_{\mathcal{H}}<\varepsilon,$\thinspace\ and then a
number $\,N\geq1$\thinspace\ such that $\,\Vert\varphi\,\mathsf{1}%
_{K}\,\mathsf{1}_{\{\varphi>N\}}\Vert_{\mathcal{H}}<\varepsilon.$%
\thinspace\ Then%
\begin{align}
\Vert S_{t}\varphi-\varphi\Vert_{\mathcal{H}}\;  &  \leq\;\big\|S_{t}%
(\varphi\,\mathsf{1}_{K^{\mathrm{c}}})-\varphi\,\mathsf{1}_{K^{\mathrm{c}}%
}\big\|_{\mathcal{H}}+\big\|S_{t}(\varphi\,\mathsf{1}_{K}\,\mathsf{1}%
_{\{\varphi>N\}})-\varphi\,\mathsf{1}_{K}\,\mathsf{1}_{\{\varphi
>N\}}\big \|_{\mathcal{H}}\nonumber\\[1pt]
&  \qquad\qquad\qquad+\;\big\|S_{t}(\varphi\,\mathsf{1}_{K}\,\mathsf{1}%
_{\{\varphi\leq N\}})-\varphi\,\mathsf{1}_{K}\,\mathsf{1}_{\{\varphi\leq
N\}}\big\|_{\mathcal{H}}\nonumber\\[1pt]
&  \leq\;C\,\varepsilon+\big\|S_{t}(\varphi\,\mathsf{1}_{K}\,\mathsf{1}%
_{\{\varphi\leq N\}})-\varphi\,\mathsf{1}_{K}\,\mathsf{1}_{\{\varphi\leq
N\}}\big\|_{\mathcal{H\,}},
\end{align}
where in the last step we used twice Lemma \ref{new-mink} with $\,\beta
=0.$\thinspace\ Thus, for the rest of the proof we may assume that $\,\varphi
$\thinspace\ is bounded by $\,N\geq1$\thinspace\ and vanishes outside a
compact set $\,K\subset\mathsf{R}^{d}.$\thinspace\ That is, from now on in
this proof we assume that $\,\mathsf{R}^{d}$\thinspace\ is replaced by
$\,K$\thinspace\ in the definition of $\,\mathcal{H}.$\medskip

\noindent2$^{\circ}$ (\emph{Passing to a continuous function}). Fix
$\,\varepsilon\in(0,1].$\thinspace\ Choose a continuous non-negative function
$\,f_{\varepsilon}\leq N$\thinspace\ (on $\,K)$\thinspace\ such that
$\,\varphi=f_{\varepsilon}$\thinspace\ on a measurable set $\,A_{\varepsilon
}\subseteq K$\thinspace\ satisfying $\,\ell(A_{\varepsilon}^{\mathrm{c}}%
)\leq\varepsilon.$\thinspace\ Then, again by twice applying Lemma
\ref{new-mink} with $\,\beta=0,$%
\begin{gather}
\Vert S_{t}\varphi-\varphi\Vert_{\mathcal{H}}\;\leq\;\big\|S_{t}%
(\varphi\,\mathsf{1}_{A_{\varepsilon}^{\mathrm{c}}})-\varphi\,\mathsf{1}%
_{A_{\varepsilon}^{\mathrm{c}}}\big\|_{\mathcal{H}}+\big\|S_{t}(f_{\varepsilon
}\,\mathsf{1}_{A_{\varepsilon}})-f_{\varepsilon}\,\mathsf{1}_{A_{\varepsilon}%
}\big\|_{\mathcal{H}}\nonumber\\[1pt]
\leq\;C\,\Vert\varphi\,\mathsf{1}_{A_{\varepsilon}^{\mathrm{c}}}%
\Vert_{\mathcal{H}}+C\,\Vert f_{\varepsilon}\,\mathsf{1}_{A_{\varepsilon
}^{\mathrm{c}}}\Vert_{\mathcal{H}}+\Vert S_{t}f_{\varepsilon}-f_{\varepsilon
}\Vert_{\mathcal{H\,}}. \label{seco}%
\end{gather}
For $\,x\in K$\thinspace\ fixed, $\,S_{t}f_{\varepsilon}(x)\rightarrow
f_{\varepsilon}(x)$\thinspace\ as $\,t\downarrow0,$\thinspace\ \vspace{1pt}%
\begin{equation}
\sup_{t\geq0}\Vert S_{t}f_{\varepsilon}\Vert_{\infty}\;\leq\;\Vert
f_{\varepsilon}\Vert_{\infty}\;<\;\infty\label{sup.S}%
\end{equation}
(with $\,\Vert\cdot\Vert_{\infty}$\thinspace\ denoting the supremum norm), and
$\,\phi$\thinspace\ is integrable on $\,K.$\thinspace\ Hence, by dominated
convergence, the third term in (\ref{seco}) will vanish as $\,t\downarrow
0,$\thinspace\ for fixed $\,\varepsilon.$\thinspace\ On the other hand,
$\,\Vert\varphi\,\mathsf{1}_{A_{\varepsilon}^{\mathrm{c}}}\Vert_{\mathcal{H}}%
$\thinspace\ converges to $0$ as $\,\varepsilon\downarrow0.$\thinspace
\ Finally, the same is true for $\,\Vert f_{\varepsilon}\,\mathsf{1}%
_{A_{\varepsilon}^{\mathrm{c}}}\Vert_{\mathcal{H}}$\thinspace\ since
$\,f_{\varepsilon}\leq N.$\thinspace\ This completes the proof.\hfill
\rule{0.5em}{0.5em}

\subsection{The fundamental solutions $P^{\alpha}$\label{SS.P.alpha}}

Fix $\,\alpha\in\mathsf{R}.$\thinspace\ We now introduce the fundamental
solutions $P^{\alpha}=P^{\alpha}(t;x,y)$ of the heat equation with one-point
potential $\,\delta_{0\,}^{(\alpha)},$\thinspace\ that is of equation
(\ref{singular'}).\medskip

\noindent1$^{\circ}$ ($d=3$). Based on \cite[formula array (3.4)]%
{AlbeverioBrzezniakDabrowski1995}, for $d=3,$ we can \emph{define}%
\begin{align}
P^{\alpha}(t;x,y)\;:=\; &  \;P(t;x,y)\,+\,\frac{2t}{|x|\,|y|}\,P\!\left(
t;|x|+|y|_{\!_{\!_{\,}}}\right)  \label{p-alpha-3}\\[2pt]
&  -\;\frac{8\pi\alpha\,t}{|x|\,|y|}\int_{0}^{\infty}\mathrm{d}u\;\mathrm{e}%
^{-4\pi\alpha u}P\!\left(  t;u+|x|+|y|_{\!_{\!_{\,}}}\right)  \!,\nonumber
\end{align}
$t>0,$\thinspace\ $x,y\neq0,$\thinspace\ where with a slight abuse of notation
for heat kernel $\,P$,%
\begin{equation}
P(t;r)\;:=\;(4\pi t)^{-d/2}\exp\!\left(  -r^{2}/4t\right)  \!,\qquad
t,r>0.\label{abuse}%
\end{equation}
Note that the term in (\ref{p-alpha-3}) involving the integral is always
finite and that it disappears for $\,\alpha=0.$ Otherwise, using the
substitution $\,|\alpha|u\rightarrow u$\thinspace\ (for $\,\alpha\neq
0)$\thinspace\ one realizes that $\,P^{\alpha}(t;x,y)$\thinspace\ is
continuous and decreasing in $\alpha,$ and that $\,P^{\alpha}\downarrow
P,$\thinspace\ the heat kernel, pointwise as $\alpha\uparrow\infty,$%
\thinspace\ whereas $\,P^{\alpha}\uparrow\infty$\thinspace\ pointwise as
$\,\alpha\downarrow-\infty.$\medskip

\noindent2$^{\circ}$ ($d=2$). On the other hand, by \cite[formula
(3.15)]{AlbeverioBrzezniakDabrowski1995}, for $d=2$, we may \emph{define}%
\begin{align}
&  P^{\alpha}(t;x,y)\;:=\;P(t;x,y)\,+\,\frac{\sqrt{4\pi t}}{\sqrt{|x|\,|y|}%
}\,P\!\left(  t;|x|+|y|_{\!_{\!_{\,}}}\right)  \times\label{p-alpha-2}\\[2pt]
&  \quad\int_{0}^{\infty}\mathrm{d}u\;\frac{t^{u}\,\mathrm{e}^{-\alpha u}%
}{\Gamma(u)}\int_{0}^{\infty}\mathrm{d}r\;\frac{r^{u-1}\,\mathrm{e}%
^{-(|x|+|y|)^{2}/4tr}}{(r+1)^{u+1/2}}\,\tilde{K}_{0}\Big(\frac{|x|\,|y|}%
{2t}(r+1)\!\Big)\nonumber
\end{align}
$t>0,$\thinspace\ $x,y\neq0,$\thinspace\ where $\Gamma$ is the Gamma
function,
\begin{equation}
\Gamma(u)\;:=\;\int_{0}^{\infty}\mathrm{d}s\;s^{u-1}\mathrm{e}^{-s},\qquad
u>0, \label{Gamma}%
\end{equation}
and
\begin{equation}
\tilde{K}_{0}(z)\;:=\;\mathrm{e}^{z}\left(  2z/\pi\right)  ^{1/2}%
K_{0}(z),\qquad z\geq0, \label{McD}%
\end{equation}
with $K_{0}\geq0$ the Macdonald function of order $0$. In other words, $K_{0}$
is the modified Bessel function of the third kind, of order $0$. See
\cite[p.109]{Lebedev1965}.

Recall that $\,P^{a}$\thinspace\ $(d=2,3)$\thinspace\ is the family of
fundamental solutions to equation (\ref{singular'}), computed in
\cite{AlbeverioBrzezniakDabrowski1995}. Since $\Delta^{\!(\alpha)}$ is a
self-adjoint extension of $\Delta$\thinspace\ on $\,\mathsf{\dot{R}}^{d}%
,$\thinspace\ the kernel $\,P^{\alpha}$\thinspace\ solves the heat equation on
$\,(0,\infty)\times\mathsf{\dot{R}}^{d}:$

\begin{cor}
[\textbf{Solutions of the heat equation}]\label{Cor.heat}Let $\,d=2,3$%
\thinspace\ and $\,\alpha\in\mathsf{R}.$\thinspace\ Then%
\begin{equation}
\frac{\partial}{\partial t}P^{\alpha}(t;x,y)\ =\ \Delta P^{\alpha}%
(t;x,y)\quad\text{on}\quad(0,\infty)\times\mathsf{\dot{R}}^{d},
\end{equation}
where the Laplacian acts on $\,x$\thinspace\ $($or $\,y,$\thinspace
\ respectively). In particular, $\,(t,x,y)\mapsto P^{\alpha}(t;x,y)$%
\thinspace\ is jointly continuous on $\,(0,\infty)\times\mathsf{\dot{R}}^{d}.$
\end{cor}

Let $S^{\alpha}=\{S_{t}^{\alpha}:\,t\geq0\}$ denote the semigroup
corresponding to the kernel $P^{\alpha},$\thinspace\ $\alpha\in\mathsf{R}.$

\subsection{Bounds of $P^{\alpha}$}

In this subsection we will derive some bounds for the kernels $\,P^{\alpha}%
$\thinspace\ introduced in (\ref{p-alpha-3}) and (\ref{p-alpha-2}),
respectively. To this end, we set%
\begin{equation}
\bar{P}(t;x,y)\;:=\;t^{-1/2}\,\phi(x)\,\phi(y)\,\mathrm{e}^{-|x|^{2}%
/4t}\,\,\mathrm{e}^{-|y|^{2}/4t}\,{,} \label{bardef}%
\end{equation}
for $\,t>0$\thinspace\ and $\,x,y\neq0$\thinspace\ [recall the weight and
reference function $\,\phi$\thinspace\ from (\ref{ref.fct})].

\begin{lemma}
[$P^{\alpha}$ bound]\label{L.kern.b}Let $\,d=2,3.$\thinspace\ For each
$\alpha\in\mathsf{R}$\thinspace\ and $\,T>0,$\thinspace\ there is a constant
$\,C_{\mathrm{\ref{L.kern.b}}}=C_{\mathrm{\ref{L.kern.b}}}(d,\alpha
,T)$\thinspace\ such that%
\begin{equation}
P(t;x,y)\;\leq\;P^{\alpha}(t;x,y)\;\leq\;P(t;x,y)+C_{\mathrm{\ref{L.kern.b}}%
}\,\bar{P}(t;x,y) \label{p-diff}%
\end{equation}
for all $\,t\in(0,T]$\thinspace\ and $\,x,y\neq0.$
\end{lemma}

%

\proof
1$^{\circ}$ ($d=3$). By the arguments after (\ref{abuse}), for $\alpha\geq0,$%
\begin{equation}
P\;\leq\;P^{\alpha}\;\leq\;P^{0}\;\leq\;P+2\,(4\pi)^{-3/2}\,\bar{P},
\label{p-diff.0}%
\end{equation}
since $\,$%
\begin{equation}
\left(  |x|+|y|_{\!_{\!_{\,}}}\right)  ^{2}\;\geq\;|x|^{2}+|y|^{2}%
.\,\ \label{bin}%
\end{equation}
So we will restrict our attention to $\,\alpha<0.$ Abbreviating%
\begin{equation}
-4\pi\alpha\;=:\;\frac{r}{2}\;>\;0\quad\text{and}\quad|x|+|y|\;=:\;R\;\geq\;0
\label{set}%
\end{equation}
and using the last inequality in (\ref{p-diff.0}) and (\ref{bin}), it suffices
to verify that%
\begin{equation}
\int_{0}^{\infty}\mathrm{d}u\;\,\frac{r}{2}\,\mathrm{e}^{\frac{r}{2}%
u}\,P(t;u+R)\;\leq\;C_{\mathrm{(\ref{verify})}}\,P(t;R) \label{verify}%
\end{equation}
[recall notation (\ref{abuse})] with a positive constant
$\,C_{\mathrm{(\ref{verify})}}=C_{\mathrm{(\ref{verify})}}(T,r)$%
\thinspace\ independent of $\,t$\thinspace\ and $\,R.$\thinspace\ Fix any
\begin{equation}
u_{0}\;>\;rT\quad\text{and put}\quad u_{1}\;:=\;u_{0}-rT\;>\;0. \label{u0u1}%
\end{equation}
Consider first the integral in (\ref{verify}) restricted to $\,u\in
\lbrack0,u_{0}].$\thinspace\ Here we can use $\,P(t;u+R)\leq P(t;R)$%
\thinspace\ and the fact that%
\begin{equation}
\int_{0}^{u_{0}}\mathrm{d}u\;\frac{r}{2}\,\mathrm{e}^{\frac{r}{2}u}%
\;\leq\;\mathrm{e}^{\frac{r}{2}u_{0}},
\end{equation}
resulting in a positive constant independent of $\,t$\thinspace\ and $\,R.$ It
remains to deal with%
\begin{subequations}
\begin{align}
&  \int_{u_{0}}^{\infty}\mathrm{d}u\;\frac{r}{2}\,\mathrm{e}^{\frac{r}{2}%
u}\,(4\pi t)^{-3/2}\,\,\mathrm{e}^{-(u+R)^{2}/4t}\,\label{integral}\\[2pt]
&  =\;\frac{r}{2}\,(4\pi t)^{-3/2}\,\,\mathrm{e}^{-(2rRt-r^{2}t^{2})/4t}%
\int_{u_{0}}^{\infty}\!\mathrm{d}u\;\,\mathrm{e}^{-(u+R-rt)^{2}/4t}
\label{integral'}%
\end{align}
for $\,0<t\leq T.$\thinspace\ The exponential factor in front of the integral
in (\ref{integral'}) is bounded by%
\end{subequations}
\begin{equation}
\mathrm{e}^{r^{2}T/4}\;=:\;C_{(\mathrm{\ref{const}})} \label{const}%
\end{equation}
which is a positive constant independent of $\,t$ and $R.\,$ Substituting
$\,u+R-rt\rightarrow u$\thinspace\ and recalling notation (\ref{u0u1}), the
integral in (\ref{integral'}) can be bounded by%
\begin{equation}
\int_{R+u_{1}}^{\infty}\!\mathrm{d}u\;\,\mathrm{e}^{-u^{2}/4t}\;\leq
\;\frac{2T}{u_{1}}\int_{R}^{\infty}\!\mathrm{d}u\;\frac{u}{2t}\,\mathrm{e}%
^{-u^{2}/4t}\;=\;\frac{2T}{u_{1}}\,\,\mathrm{e}^{-R^{2}/4t}.
\end{equation}
Thus for the integral in (\ref{integral}) we found the bound%
\begin{equation}
C_{(\mathrm{\ref{const}})}\,\frac{r}{2}\,\frac{2T}{u_{1}}\,P(t;R),
\end{equation}
which finishes the proof in the case $\,d=3.$\medskip

\noindent2$^{\circ}$ ($d=2$). Recall definition (\ref{McD}) of $\,\tilde
{K}_{0\,}.$\thinspace\ According to \cite[after (3.14)]%
{AlbeverioBrzezniakDabrowski1995},
\begin{equation}
\lim_{z\rightarrow\infty}\tilde{K}_{0}(z)\;=\;1. \label{toi}%
\end{equation}
Consulting \cite[Section 1.15, equation (1.66)]{Tranter1969}, we find that%
\begin{equation}
K_{0}(z)\;\sim\;-\gamma-\log(z/2)\ \sim\ -\log z\quad\text{as }\,z\downarrow0,
\end{equation}
where $\,\gamma$\thinspace\ is Euler's constant. Therefore,
\begin{equation}
\lim_{z\downarrow0}\tilde{K}_{0}(z)\;=\;\lim_{z\downarrow0}\left[
\mathrm{e}^{z}(2z/\pi)^{1/2}\,\log z_{\!_{\!_{\,}}}\right]  \,=\;0.
\label{to0}%
\end{equation}
Since $\,\tilde{K}_{0}$\thinspace\ is continuous, relations (\ref{toi}) and
(\ref{to0}) together give%
\begin{equation}
\Vert\tilde{K}_{0}\Vert_{\infty}\,<\,\infty.
\end{equation}

Fix $\,\alpha\in\mathsf{R}$\thinspace\ and consider $\,0<t\leq T.$%
\thinspace\ We may assume that $\,T\geq1.$\thinspace\ We start by estimating
the inner integral appearing on the right hand side of definition
(\ref{p-alpha-2}) of $\,P^{\alpha}$. For $\,u>0,$
\begin{gather}
\int_{0}^{\infty}\mathrm{d}r\;\frac{r^{u-1}\mathrm{e}^{-(|x|+|y|)^{2}/4tr}%
}{(r+1)^{u+1/2}}\,\tilde{K}_{0}\left(  \frac{|x|\,|y|}{2t}(r+1)\!\right)
\label{e-i-1}\\
\leq\;\Vert\tilde{K}_{0}\Vert_{\infty}\int_{0}^{\infty}\mathrm{d}%
r\;\frac{r^{u-1}}{(r+1)^{u+1/2}}\,.\nonumber
\end{gather}
If $\,r\geq1,$ drop the $1$ in the denominator, otherwise drop the
$\,r$\thinspace\ there. Thus, for the inner integral in (\ref{p-alpha-2}) we
find the bound%
\begin{equation}
\Vert\tilde{K}_{0}\Vert_{\infty}\,\left[  2+1/u\right]  \!.\label{e-i-2}%
\end{equation}
Using this bound, we turn to the outer integral of (\ref{p-alpha-2}). For the
Gamma function $\Gamma$ of (\ref{Gamma}), Stirling's formula gives
\begin{equation}
\Gamma(u)\;\sim\;\sqrt{2\pi}\,(u-1)^{u-1/2}\,\mathrm{e}^{-u+1}\quad\text{as
}\,u\uparrow\infty.
\end{equation}
It follows that, for some constant $\,C_{(\mathrm{\ref{e-i-3}})}%
=C_{(\mathrm{\ref{e-i-3}})}(T,\alpha),$%
\begin{equation}
\int_{1}^{\infty}\mathrm{d}u\;\text{\thinspace}\frac{t^{u}\,\mathrm{e}%
^{-\alpha u}}{\Gamma(u)}\;\leq\;C_{(\mathrm{\ref{e-i-3}})},\qquad0\leq t\leq
T.\label{e-i-3}%
\end{equation}
Next, using integration by parts, we estimate $\,u\,\Gamma(u)$\thinspace\ for
$\,u\in(0,1]:$%
\begin{equation}
u\,\Gamma(u)\;=\;\int_{0}^{\infty}\mathrm{d}s\;us^{u-1}\mathrm{e}%
^{-s}\;=\;\int_{0}^{\infty}\mathrm{d}s\;s^{u}\,\mathrm{e}^{-s}\;\geq
\;\mathrm{e}^{-1}\int_{0}^{1}\mathrm{d}s\;s\;=:\;C_{(\mathrm{\ref{uG}}%
)}.\label{uG}%
\end{equation}
Finally, for some constant $\,C_{(\mathrm{\ref{C01}})}=C_{(\mathrm{\ref{C01}%
})}(T,\alpha),$\thinspace\ since $\,T\geq1,$
\begin{equation}
\int_{0}^{1}\mathrm{d}u\;\,\frac{t^{u}\,\mathrm{e}^{-\alpha u}}{u\,\Gamma
(u)}\;\leq\;\frac{1}{C_{(\mathrm{\ref{uG}})}}T\,\mathrm{e}^{|\alpha
|}\;=:\;C_{(\mathrm{\ref{C01}})}.\label{C01}%
\end{equation}
Altogether, we found that the double integral appearing on the right hand side
of definition (\ref{p-alpha-2}) of $\,P^{\alpha}$\thinspace\ is bounded by a
constant depending only on $\,\alpha,T$. This gives estimate (\ref{p-diff})
also in the case $\,d=2,$\thinspace\ since $\,\tilde{K}_{0}\geq0,$%
\thinspace\ finishing the proof of Lemma \ref{L.kern.b}.\hfill
\rule{0.5em}{0.5em}

\subsection{Strong continuity of $S^{\alpha}$}

We abbreviate
\begin{equation}
\bar{S}_{t}\varphi\,(x)\;:=\;\int_{\mathsf{R}^{d}}\mathrm{d}y\;\varphi
(y)\,\bar{P}(t;x,y),\qquad t>0,\quad x\neq0, \label{barSdef}%
\end{equation}
with $\bar{P}$ from (\ref{bardef}), as long as the right hand side expression
makes sense. The estimates (\ref{p-diff}) and Minkowski's inequality then
imply that
\begin{equation}
\Vert S_{t}\varphi\Vert_{\mathcal{H}}\;\leq\;\Vert S_{t}^{\alpha}\varphi
\Vert_{\mathcal{H}}\;\leq\;\Vert S_{t}\varphi\Vert_{\mathcal{H}}%
+C_{\mathrm{\ref{L.kern.b}}}\,\Vert\bar{S}_{t}\varphi\Vert_{\mathcal{H\,}%
},\qquad0<t\leq T, \label{mink-1}%
\end{equation}
for those $\varphi$ for which the right hand side of (\ref{mink-1}) is
meaningful and finite.

\begin{lemma}
[Estimate of $\bar{S}$ in case of an additional singularity]%
\label{new-mink.II}\hspace{-2.1pt}Let $\,d$ $=2,3$\thinspace\ as well as
$\,0\leq\beta\leq1,$\thinspace\ and assume%
\begin{equation}
\frac{1}{1-\beta(d-1)/(d+1)}\;<\;\varrho\;<\;\frac{d+1}{d-1}\,. \label{back}%
\end{equation}
Then there is a constant $\,C_{\mathrm{\ref{new-mink.II}}\,}%
=C_{\mathrm{\ref{new-mink.II}}\,}(d,\beta,\varrho)$\thinspace\ such that for
all $\,\varphi\in\mathcal{H}=\mathcal{H}^{\varrho},$%
\begin{equation}
\left\Vert \bar{S}_{t}(\varphi\phi^{\beta})\right\Vert _{\mathcal{H}}%
^{\varrho}\;\leq\;C_{\mathrm{\ref{new-mink.II}}}\,\varepsilon(t,\varphi
)\,t^{-\beta\varrho(d-1)/4}\,\Vert\varphi\Vert_{\mathcal{H}}^{\varrho
}\,,\qquad t>0, \label{bar-norm-est}%
\end{equation}
where $\,0\leq\varepsilon(t,\varphi)\leq1\,\ $and\/ $\,\varepsilon
(t,\varphi)\rightarrow0$\thinspace\ as $\,t\downarrow0.$\smallskip
\end{lemma}

\begin{remark}
[Restriction to infinite variance branching if $d=3$]\label{R.beta}\hfill Note
\newline that in dimension $\,d=3$\thinspace\ condition (\ref{back}) can only
be satisfied for some $\,\varrho$\thinspace\ if $\,\beta<1$\thinspace
\ holds.\hfill$\Diamond$\smallskip
\end{remark}

\noindent\emph{Proof of Lemma}\/ \ref{new-mink.II}. \thinspace This time we
work with the measures%
\begin{equation}
\mu_{t}(\mathrm{d}y)\;:=\;t^{-\kappa}\,\mathrm{e}^{-|y|^{2}/4t}\,\phi
^{\lambda}(y)\,\mathrm{d}y,\qquad t>0, \label{def.mu.t}%
\end{equation}
on $\,\mathsf{R}^{d},$\thinspace\ where%
\begin{equation}
\kappa\;:=\;\frac{d+1}{4}-\frac{(d-1)\beta\varrho}{4(\varrho-1)}%
\;>0\quad\text{and}\quad\lambda\;:=\;1+\beta\varrho/(\varrho-1).
\end{equation}
Note that the measures $\mu_{t}$ have a $t$--independent total mass%
\begin{equation}
\Vert\mu_{t}\Vert\;=\;\int_{\mathsf{R}^{d}}\mathrm{d}y\;\mathrm{e}%
^{-|y|^{2}/4}\,\phi^{\lambda}(y)\;=:\;C_{(\text{\textrm{\ref{CNr}}}%
)}=C_{(\text{\textrm{\ref{CNr}}})}(d,\beta,\varrho), \label{CNr}%
\end{equation}
which is finite by the left hand inequality in assumption (\ref{back}). Then,
by our definition (\ref{bardef}) of $\,\bar{P},$\thinspace\ for $\,t>0$%
\thinspace\ and $\,x\neq0,$%
\begin{align}
&  \left\vert \bar{S}_{t}(\varphi\phi^{\beta})\,(x)\right\vert ^{\varrho
}\label{mink-4.5}\\
&  =\;t^{-\varrho/2+\kappa\varrho}\,\phi^{\varrho}(x)\,\mathrm{e}%
^{-\varrho|x|^{2}/4t}\,\bigg|\int_{\mathsf{R}^{d}}\mu_{t}(\mathrm{d}%
y)\,\phi^{-\lambda+\beta+1}(y)\,\varphi(y)\bigg|^{\varrho}.\nonumber
\end{align}
By (\ref{Hoelder}) and the definition of $\,\mu_{t}$\thinspace\ we may
continue with%
\begin{align}
&  \left\vert \bar{S}_{t}(\varphi\phi^{\beta})(x)\right\vert ^{\varrho}\\
&  \leq\;t^{-\varrho/2+\kappa\varrho}\,\phi^{\varrho}(x)\,\mathrm{e}%
^{-\varrho|x|^{2}/4t}\,C_{(\text{\textrm{\ref{CNr}}})}^{\varrho-1}%
\,t^{-\kappa}\int_{\mathsf{R}^{d}}\!\mathrm{d}y\;\mathrm{e}^{-|y|^{2}%
/4t}\,\phi(y)\,|\varphi|^{\varrho}(y),\nonumber
\end{align}
since $\,(-\lambda+\beta+1)\varrho+\lambda=1.$\thinspace\ We may assume that
$\,\varphi\neq0.$\thinspace\ Define%
\begin{equation}
\varepsilon(t,\varphi)\;:=\;\frac{1}{\Vert\varphi\Vert_{\mathcal{H}}^{\varrho
}}\,\int_{\mathsf{R}^{d}}\mathrm{d}y\;\phi(y)\,\mathrm{e}^{-|y|^{2}%
/4t}\,|\varphi|^{\varrho}(y).
\end{equation}
Note that $\,0<\varepsilon(t,\varphi)\leq1$\thinspace\ and that $\,\varepsilon
(t,\varphi)\rightarrow0$\thinspace\ as $\,t\downarrow0,$\thinspace\ by
dominated convergence. Consequently,%
\begin{equation}
\left\vert \bar{S}_{t}(\varphi\phi^{\beta})(x)\right\vert ^{\varrho}%
\,\leq\;t^{-\varrho/2+\kappa\varrho}\,\phi^{\varrho}(x)\,\mathrm{e}%
^{-\varrho|x|^{2}/4t}\,C_{(\text{\textrm{\ref{CNr}}})}^{\varrho-1}%
\,t^{-\kappa}\,\varepsilon(t,\varphi)\,\Vert\varphi\Vert_{\mathcal{H\,}%
}^{\varrho}. \label{conse}%
\end{equation}
Therefore,%
\[
\left\Vert \bar{S}_{t}(\varphi\phi^{\beta})\right\Vert _{\mathcal{H}}%
^{\varrho}\;\leq\;C_{(\text{\textrm{\ref{CNr}}})}^{\varrho-1}\,\varepsilon
(t,\varphi)\,t^{-\varrho/2+\kappa\varrho-\kappa}\,\Vert\varphi\Vert
_{\mathcal{H}}^{\varrho}\int_{\mathsf{R}^{d}}\mathrm{d}x\;\phi^{\varrho
+1}(x)\,\mathrm{e}^{-\varrho|x|^{2}/4t}.
\]
But the latter integral is finite since $\,-(\varrho+1)(d-1)/2+d>0$%
\thinspace\ by the right hand inequality in assumption (\ref{back}). Moreover,
using a change of variables, the integral gives an additional factor
$\,t^{d/2-(\varrho+1)(d-1)/4},$\thinspace\ so that the whole $t$--term equals
$\,t^{-\beta\varrho(d-1)/4}.$\thinspace\ This finishes the proof.\hfill
\rule{0.5em}{0.5em}\bigskip

Since condition (\ref{back}) is stronger than (\ref{ass}), combining Lemmas
\ref{new-mink} and \ref{new-mink.II} with inequality (\ref{mink-1}) gives the
following result.

\begin{cor}
[Estimate of $S^{\alpha}$ in case of an additional singularity]%
\label{Cor.est.sing}\quad\newline Let $\,d=2,3$\thinspace\ as well as
$\,0\leq\beta\leq1.$\thinspace\ Suppose\/ \emph{(\ref{back}). }To each
$\,T>0$\thinspace\ there is a constant $\,C_{\mathrm{\ref{Cor.est.sing}}%
\,}=C_{\mathrm{\ref{Cor.est.sing}}\,}(d,T,\alpha,\beta,\varrho)$%
\thinspace\ such that for all $\,\varphi\in\mathcal{H}=\mathcal{H}^{\varrho},$%
\begin{equation}
\left\Vert S_{t}^{\alpha}(\varphi\phi^{\beta})\right\Vert _{\mathcal{H}}%
\;\leq\;C_{\mathrm{\ref{Cor.est.sing}}}\,t^{-\beta(d-1)/4}\,\Vert\varphi
\Vert_{\mathcal{H}}\,,\qquad0<t\leq T.\label{est.sing}%
\end{equation}
In particular,%
\begin{equation}
\sup_{t\leq T}\Vert S_{t}^{\alpha}\varphi\Vert_{\mathcal{H}}\;\leq
\;C_{\mathrm{\ref{Cor.est.sing}}}\,\Vert\varphi\Vert_{\mathcal{H}}%
\;<\;\infty,\qquad\varphi\in\mathcal{H}.\label{bounded}%
\end{equation}

\end{cor}

Here is another consequence of Lemma \ref{new-mink.II}:

\begin{cor}
[\textbf{Strong continuity of }$S^{\alpha}$]\label{C.scSalpha}Let
$\,d=2,3.$\thinspace\ For each $\,\alpha\in\mathsf{R},$\thinspace\ the
semigroup\thinspace\ $S^{\alpha}$\thinspace\ acting on $\,\mathcal{H=H}%
^{\varrho}$\thinspace\ with $\,\varrho\in\left(  1,\,(d+1)/(d-1)_{\!_{\!_{\,}%
}}\right)  $\thinspace\ is strongly continuous.
\end{cor}

%

\proof
Fix $\,\varphi\in\mathcal{H}.$\thinspace\ By linearity, we may additionally
assume that $\,\varphi\geq0.$\thinspace\ Consider $\,0<t\leq T.$%
\thinspace\ Decompose
\begin{equation}
S_{t}^{\alpha}\varphi\;=\;S_{t}\varphi+(S_{t}^{\alpha}-S_{t})\varphi,
\label{dec}%
\end{equation}
where by Lemma \ref{L.kern.b},%
\begin{equation}
0\;\leq\;(S_{t}^{\alpha}-S_{t})\varphi\;\;\leq\;C_{\mathrm{\ref{L.kern.b}}%
}\,\bar{S}_{t}\varphi,
\end{equation}
implying%
\begin{equation}
\left\|  (S_{t}^{\alpha}-S_{t})\varphi_{\!_{\!_{\,}}}\right\|  _{\mathcal{H}%
}\;\;\leq\;C_{\mathrm{\ref{L.kern.b}}}\,\Vert\bar{S}_{t}\varphi\Vert
_{\mathcal{H\,}}. \label{bound}%
\end{equation}
\textrm{F}rom (\ref{dec}) and (\ref{bound})%
\begin{gather}
\Vert S_{t}^{\alpha}\varphi-\varphi\Vert_{\mathcal{H}}\;\leq\;\Vert
S_{t}\varphi-\varphi\Vert_{\mathcal{H}}+\left\|  (S_{t}^{\alpha}-S_{t}%
)\varphi_{\!_{\!_{\,}}}\right\|  _{\mathcal{H}}\nonumber\\[1pt]
\leq\;\Vert S_{t}\varphi-\varphi\Vert_{\mathcal{H}}+C_{\mathrm{\ref{L.kern.b}%
}}\,\Vert\bar{S}_{t}\varphi\Vert_{\mathcal{H\,}}. \label{sec}%
\end{gather}
But by Lemma \ref{new-mink.II} with $\,\beta=0,$\thinspace\ the second term in
(\ref{sec}) goes to $0$ as $\,t\downarrow0,$\thinspace\ whereas the first term
does by Lemma \ref{L.cont.heat}. By (\ref{bounded}), this finishes the
proof.\hfill\rule{0.5em}{0.5em}

\subsection{$S^{\alpha}$ as a flow on $\Phi$\label{SS.SalphaPhi}}

Recall our set $\,\Phi$\thinspace\ of continuous non-negative test functions
introduced in Subsection \ref{SS.Prel}. \textrm{F}rom the proof of Lemma
\ref{new-mink.II} we also get the following result:

\begin{cor}
[$S^{\alpha}$ \textbf{bound}]\label{Cor.flow.bound}Let $\,d=2,3,$%
\thinspace\ assume $\,\varrho\in\left(  1,\,(d+1)/(d-1)_{\!_{\!_{\,}}}\right)
\!,$\thinspace\ and that $\,\varphi\in\mathcal{H}^{\varrho}$\thinspace
\ satisfies\thinspace\footnote{$^{)}$ Of course, an inequality on an element
$\,\varphi\in\mathcal{H}$\thinspace\ means that the inequality holds for each
representative in Lebesgue almost all points.}\/ \emph{(\ref{phi.bound}).}
Then, to each $\,T>0,$\thinspace\ there is a constant
$\,C_{\mathrm{\ref{Cor.flow.bound}}}=C_{\mathrm{\ref{Cor.flow.bound}}%
}(d,T,\alpha,\varrho,\varphi)$\thinspace\ such that%
\begin{equation}
0\;\leq\;S_{t}^{\alpha}\varphi\;\leq\;C_{\mathrm{\ref{Cor.flow.bound}}%
}\,\big (1+t^{-1/2+(d+1)(\varrho-1)/4\varrho}\big)\phi,\qquad0<t\leq T.
\label{flow.bound}%
\end{equation}
In particular, $\,S_{t}^{\alpha}\varphi\in\Phi,$\thinspace\ for all $\,t>0.$
\end{cor}

%

\proof
\textrm{F}rom Lemma \ref{L.kern.b},%
\begin{equation}
0\;\leq\;S_{t}^{\alpha}\varphi\;\leq\;S_{t}\varphi+C_{\mathrm{\ref{L.kern.b}}%
}\,\bar{S}_{t}\varphi. \label{from}%
\end{equation}
Moreover, by assumption (\ref{phi.bound}) on $\,\varphi$\thinspace\ and by
Lemma \ref{L.hf.est},%
\begin{equation}
S_{t}\varphi\;\leq\;C_{\mathrm{(\ref{phi.bound})}}S_{t}\phi\;\leq\;C\phi.
\end{equation}
On the other hand, raising estimate (\ref{conse}) (with $\,\beta=0$%
\thinspace\ there, implying $\,\kappa=(d+1)/4$\thinspace\ and $\,\lambda
=1)$\thinspace\ into the power $\,1/\varrho$\thinspace\ gives%
\begin{equation}
\bar{S}_{t}\varphi\;\leq\;C\,t^{-1/2\,+\,(d+1)(\varrho-1)/4\varrho}%
\,\phi\,\Vert\varphi\Vert_{\mathcal{H\,}}. \label{raising}%
\end{equation}
Putting together (\ref{from})\thinspace--\thinspace(\ref{raising}) yields
(\ref{flow.bound}). Finally, $\,(t,x)\mapsto S_{t}^{\alpha}\varphi$%
\thinspace\ is continuous on $\,(0,\infty)\times\mathsf{\dot{R}}^{d}%
,$\thinspace\ since it solves the heat equation, recall Corollary
\ref{Cor.heat}. This finishes the proof.\hfill\rule{0.5em}{0.5em}\bigskip

Combining Corollaries \ref{C.scSalpha} and \ref{Cor.flow.bound}, we get the
following result.

\begin{cor}
[$S^{\alpha}$ acting on $\Phi$]\label{Cor.on.H.bar}Let $\,d=2,3$%
\thinspace\ and $\,\varrho\in\left(  1,\,(d+1)/(d-1)_{\!_{\!_{\,}}}\right)
\!.$\thinspace\ Then $\,S^{\alpha}$\thinspace\ is a strongly continuous linear
semigroup acting on $\,\Phi=\Phi^{\varrho}.$
\end{cor}

\section{Analysis of the log-Laplace equation\label{for-both}}

The main result of this section is the well-posedness of the log-Laplace
equation (Theorem \ref{T.well}). Uniqueness follows from a contraction
argument (Lemma \ref{unique-3}). Existence is shown via a Picard iteration
(Lemmas \ref{L.cl-1}, \ref{L.cl-2}, and \ref{exist-3}), whereas non-negativity
follows using a linearized equation (Lemma \ref{from-trotter}).

\subsection{Preliminaries and purpose}

Formally, we can rewrite the log-Laplace equation (\ref{cumulant}) as the
following \emph{integral equation}\thinspace\footnote{$^{)}$ We often use
notation as $\,v(s):=v(s,\,\cdot\,).$}:
\begin{equation}
v(t,x)\;=\;S_{t}^{\alpha}\varphi\,(x)-\eta\int_{0}^{t}\mathrm{d}%
s\;S_{t-s}^{\alpha}\!\left(  v^{1+\beta}(s)\right)  (x),\label{cumulant-int}%
\end{equation}
$t\geq0,$\thinspace\ $x\neq0,$\thinspace\ (with constants $\,\alpha
\in\mathsf{R},\,\ \eta\geq0,$\thinspace\ $0<\beta\leq1,\,\ $and where
$\,\varphi\geq0\,$\ has still to be specified). \ Here in writing
$\,v^{1+\beta}$\thinspace\ we have in mind that $\,v\geq0.$\thinspace\ Note
also that this non-negativity implies the following \emph{domination:}%
\begin{equation}
0\,\leq\,v(t)\,\leq\,S_{t}^{\alpha}\varphi,\qquad t\geq0.\label{dom}%
\end{equation}
The task of this section is to verify that the log-Laplace equation
(\ref{cumulant-int}) is well-posed in $\,\Phi.$\thinspace\ 

\begin{definition}
[$\Phi$--valued solution]Let $\,\varphi\in\Phi.$\thinspace\ A measurable map
$\,t\mapsto v(t)=V_{t}\varphi$\thinspace\ of $\,\mathsf{R}_{+}$\thinspace
\ into\thinspace\ $\Phi$\thinspace\ is called a \emph{solution}\/ of
(\ref{cumulant-int}), if (\ref{cumulant-int}) is true for all $\,x\neq
0$\thinspace\ and $\,t\geq0.$\hfill$\Diamond$
\end{definition}

For convenience, we introduce the following hypothesis:

\begin{hypothesis}
[\textbf{Choice of parameters}]\label{H.choice}Let $\,\alpha\in\mathsf{R}%
,$\thinspace\ $\eta\geq0,$\thinspace\ and%
\begin{equation}
d=2,3,\quad0<\beta\leq1\quad\text{and}\quad\frac{1}{1-\beta(d-1)/(d+1)}%
\;<\;\varrho\;<\;\frac{d+1}{d-1}\,. \label{ass.rho}%
\end{equation}
Recall that for $\,d=3$\thinspace\ this requires that $\,\beta<1.$%
\hfill$\Diamond$
\end{hypothesis}

Now we are ready to state the main result of this section.

\begin{theorem}
[\textbf{Well-posedness of the log-Laplace equation}]\label{T.well}If
Hypothesis\/ \emph{\ref{H.choice}} holds, and if $\,\varphi\in\Phi,$%
\thinspace\ then equation\/ \emph{(\ref{cumulant-int})} has a unique $\Phi
$--valued solution $\,v=V\!\varphi=\left\{  V_{t}\varphi:\,t\geq0\right\}
\!.$\thinspace\ Moreover, $\,V=\{V_{t}:\,t\geq0\}$\ is a non-linear strongly
continuous semigroup acting on $\,\Phi.$
\end{theorem}

The rest of this section is devoted to the proof of this theorem.

\subsection{First properties of solutions\label{SS.first}}

Now we prepare for the uniqueness proof. Impose Hypothesis \ref{H.choice}. Fix
an integer $\,T>0$\thinspace\ for a while, and $\,\varphi\in\Phi.$%
\thinspace\ We will fix also a measurable function $\,\psi$\thinspace\ on
$\,(0,T]\times\mathsf{\dot{R}}^{d}$\thinspace\ such that%
\begin{equation}
0\;\leq\;\psi(t,x)\;\leq\;M\,(1+t^{-\kappa})\,\phi^{\beta}(x), \label{a-1}%
\end{equation}
with constants $\,M=M(T,\psi)>0$\thinspace\ and $\,$%
\begin{equation}
\kappa\;:=\;\beta/2\,-\,\beta(d+1)(\varrho-1)/4\varrho\;\in\;(0,1).
\label{kappa}%
\end{equation}

\begin{lemma}
[\textbf{Properties of the non-linear term}]\label{L.cl-0}There is a constant
$\,C_{\mathrm{\ref{L.cl-0}}}$ $=$ $C_{\mathrm{\ref{L.cl-0}}}(d,M,T,\alpha
,\beta,\varrho)$\thinspace\ such that%
\begin{equation}%
\Big\|%
\int_{0}^{t}\!\mathrm{d}s\;S_{t-s}^{\alpha}\!\left(  \psi(s)\,S_{s}^{\alpha
}\varphi\right)
\Big\|%
_{\mathcal{H}}\;\leq\ C_{\mathrm{\ref{L.cl-0}}}\,\Vert\varphi\Vert
_{\mathcal{H}}\,I(t),\qquad0<t\leq T, \label{113}%
\end{equation}
where
\begin{equation}
\infty\ >\ I(t)\ :=\ \int_{0}^{t}\!\mathrm{d}s\ (1-s^{-\kappa}%
)\,(t-s)^{-\lambda}\ \underset{t\downarrow0}{\searrow}\;0. \label{not.I.t}%
\end{equation}
Moreover, if for fixed $\,t\in(0,T],$%
\begin{equation}
N_{t}(x)\ :=\;\int_{0}^{t}\!\mathrm{d}s\;S_{t-s}^{\alpha}\!\left(
\psi(s)\,S_{s}^{\alpha}\varphi_{\!_{\!_{\,}}}\right)  (x),\qquad
x\in\mathsf{\dot{R}}^{d}, \label{finite-1}%
\end{equation}
satisfies%
\begin{equation}
N_{t}(x)\ \leq\ S_{t}^{\alpha}\varphi\,(x),\qquad x\in\mathsf{\dot{R}}^{d},
\label{ass.dom}%
\end{equation}
Then $\,N_{t}\in\Phi.$
\end{lemma}

%

\proof
First, by Corollary \ref{C.scSalpha}, we see that%
\begin{equation}
\Vert S_{s}^{\alpha}\varphi\Vert_{\mathcal{H}}\;\leq\;C\,\Vert\varphi
\Vert_{\mathcal{H\,}},\qquad0\leq s\leq T, \label{b-1}%
\end{equation}
where $\,C=C(T).$\thinspace\ Now, Corollary \ref{Cor.est.sing} states that%
\begin{equation}
\left\Vert S_{t}^{\alpha}(\varphi\,\phi^{\beta})\right\Vert _{\mathcal{H}%
}\;\leq\;C_{\mathrm{\ref{Cor.est.sing}}}\,t^{-\lambda}\,\Vert\varphi
\Vert_{\mathcal{H\,}},\qquad0<t\leq T, \label{b-2}%
\end{equation}
with $\,$%
\begin{equation}
\lambda\;:=\;\beta(d-1)/4. \label{lambda}%
\end{equation}
Applying first (\ref{b-2}) and then (\ref{b-1}), we obtain%
\begin{equation}
\left\Vert S_{t-s}^{\alpha}\!\left(  \phi^{\beta}\,S_{s}^{\alpha}%
\varphi\right)  _{\!_{\!_{\,}}}\right\Vert _{\mathcal{H}}\;\leq
\;C_{\mathrm{\ref{Cor.est.sing}}}\,(t-s)^{-\lambda}\,\Vert\varphi
\Vert_{\mathcal{H\,}},\qquad0\leq s<t\leq T. \label{b-3}%
\end{equation}
Using the fact that $\,S_{t}^{\alpha}$\thinspace\ is an integral operator with
a non-negative kernel, and exploiting assumption (\ref{a-1}), we find%
\begin{equation}
\left\Vert S_{t-s}^{\alpha}\!\left(  \psi(s)\,S_{s}^{\alpha}\varphi\right)
_{\!_{\!_{\,}}}\right\Vert _{\mathcal{H}}\;\leq\;C\,(1+s^{-\kappa
})\,(t-s)^{-\lambda}\,\Vert\varphi\Vert_{\mathcal{H\,}}. \label{112}%
\end{equation}
However, $\,I(t)$\thinspace\ from (\ref{not.I.t}) can be written as%
\begin{equation}
I(t)\;=\;\frac{t^{1-\lambda}}{1-\lambda}+t^{1-\lambda-\kappa}\int_{0}%
^{1}\!\mathrm{d}s\;s^{-\kappa}\,(1-s)^{-\lambda}. \label{scal.I}%
\end{equation}
But the positive numbers $\,\kappa\,\ $and $\,\lambda$\thinspace\ defined in
(\ref{kappa}) and (\ref{lambda}), respectively, satisfy $\,\kappa+\lambda
<1,$\thinspace\ hence (\ref{113}) and (\ref{not.I.t}) follow. Thus, the
integrals in (\ref{finite-1}) are finite for almost all $\,x.$\thinspace\ By
assumption (\ref{ass.dom}), it remains to show that $\,N_{t}(x)$%
\thinspace\ from (\ref{finite-1}) is continuous in $\,x.$

Let $\,\delta\in(0,t).$\thinspace\ Then%
\begin{equation}
\int_{0}^{t-\delta}\!\mathrm{d}s\;S_{t-s}^{\alpha}\!\left(  \psi
(s)\,S_{s}^{\alpha}\varphi\right)  (x)\ =\ S_{\delta}^{\alpha}\int
_{0}^{t-\delta}\!\mathrm{d}s\;S_{t-\delta-s}^{\alpha}\!\left(  \psi
(s)\,S_{s}^{\alpha}\varphi_{\!_{\!_{\,}}}\right)  (x).\label{part}%
\end{equation}
We already showed that the latter integral term belongs to $\,\mathcal{H}%
_{+\,}.$\thinspace\ Then by Corollary \ref{Cor.flow.bound}, the left hand side
in (\ref{part}) belongs to $\,\Phi,$\thinspace\ hence is continuous in
$\,x,$\thinspace\ for each $\,\delta.$\thinspace\ To complete the proof, it
suffices to show that%
\begin{equation}
\int_{t-\delta}^{t}\!\mathrm{d}s\;S_{t-s}^{\alpha}\!\left(  \psi
(s)\,S_{s}^{\alpha}\varphi\right)  (x)\;\underset{\delta\downarrow
0}{\longrightarrow}\;0\quad\text{uniformly in}\,\,x\in K,\label{integr}%
\end{equation}
where $\,K$\thinspace\ is any compact subset of $\,\mathsf{\dot{R}}^{d}%
,$\thinspace\ we fix from now on. Next apply Corollary~\ref{Cor.flow.bound} to
$\,S_{s}^{\alpha}\varphi$\thinspace\ together with the definition
(\ref{kappa}) of $\,\kappa$\thinspace\ to get%
\begin{equation}
0\;\leq\;S_{s}^{\alpha}\varphi\;\leq\;C_{\mathrm{\ref{Cor.flow.bound}}%
}\,(1-s^{-\kappa/\beta})\,\phi\;\leq\;C\,\phi,
\end{equation}
since $\,s$\thinspace\ in (\ref{integr}) is bounded away from $0.$ Inserting
the assumed upper bound (\ref{a-1}) on $\,\psi,$\thinspace\ for the integral
in (\ref{integr}) we find the estimate%
\begin{equation}
C\int_{t-\delta}^{t}\!\mathrm{d}s\;S_{t-s}^{\alpha}\phi^{\beta+1}%
\,(x)\ =\ C\int_{0}^{\delta}\!\mathrm{d}s\;S_{s}^{\alpha}\phi^{\beta+1}\,(x).
\end{equation}
Hence, it suffices to show that%
\begin{equation}
s\ \mapsto\ \max_{x\in K}\,S_{s}^{\alpha}\phi^{\beta+1}\,(x)\quad\text{is
integrable on }\,(0,\delta].
\end{equation}
By Lemma \ref{L.kern.b},%
\begin{equation}
S_{s}^{\alpha}\phi^{\beta+1}\ \leq\ S_{s}\phi^{\beta+1}%
+C_{\mathrm{\ref{L.kern.b}}}\,\bar{S}_{s}\phi^{\beta+1},\qquad
s>0.\label{decompose}%
\end{equation}
Now, $\,(s,x)\mapsto S_{s}\phi^{\beta+1}(x)$\thinspace\ is finite and
satisfies the heat equation on $\,[0,\delta]\times K,$\thinspace\ implying%
\begin{equation}
\sup_{(s,x)\,\in\,[0,\delta]}S_{s}\phi^{\beta+1}(x)\ <\ \infty.
\end{equation}
Turning to the second term in (\ref{decompose}), by definition (\ref{bardef}),%
\begin{equation}
\bar{S}_{s}\phi^{\beta+1}\,(x)\ =\ s^{-1/2}\,\phi(x)\,\mathrm{e}^{-|x|/4s}%
\int_{\mathsf{R}^{d}}\!\mathrm{d}y\ \mathrm{e}^{-|y|/4s}\,\phi^{\beta
+2}(y),\qquad s>0.
\end{equation}
By the substitution $\,y\mapsto y\sqrt{s},$\thinspace\ the latter integral
gives an additional power contribution to $\,s^{-1/2}.$\thinspace\ Moreover,%
\begin{equation}
\sup_{x\in K}\phi(x)\,\mathrm{e}^{-|x|/4s}\ \leq\ C\,\mathrm{e}^{-C/s},
\end{equation}
which together with $\,s^{-\lambda}$\thinspace\ is integrable on $(0,\delta],$
for each $\,\lambda.$\thinspace\ This finishes the proof.\hfill
\rule{0.5em}{0.5em}

\begin{lemma}
[\textbf{Continuity at }$t=0$]\label{eta-zero}Let $\,\varphi\in\Phi
=\Phi^{\varrho}$\thinspace\ and $\,v=V\!\varphi$\thinspace\ a $\Phi$--valued
solution to\/ \emph{(\ref{cumulant-int}).} Under Hypothesis\/
\emph{\ref{H.choice}, }for $\,T\geq0$\thinspace\ fixed, there is a constant
$\,C_{\mathrm{\ref{eta-zero}}}=C_{\mathrm{\ref{eta-zero}}}(d,T,\alpha
,\varrho)$\thinspace\ such that
\begin{equation}
\Vert V_{t}\varphi\Vert_{\mathcal{H}}\;\leq\;C_{\mathrm{\ref{eta-zero}}%
}\,\Vert\varphi\Vert_{\mathcal{H\,}},\qquad0\leq t\leq T. \label{est.eta.zero}%
\end{equation}
Moreover, $\,V\!\varphi$\thinspace\ is strongly continuous at $\,t=0,$%
\thinspace\ where $\,V_{0}\varphi=\varphi.$
\end{lemma}

\proof By domination (\ref{dom}),%
\begin{equation}
\Vert V_{t}\varphi\Vert_{\mathcal{H}}\;\leq\;\Vert S_{t}^{\alpha}\varphi
\Vert_{\mathcal{H\,}}.
\end{equation}
Now (\ref{est.eta.zero}) follows from Corollary \ref{Cor.on.H.bar}. It remains
to verify the continuity claim. Clearly, for $\,t\in(0,T],$%
\begin{equation}
\left\vert V_{t}\varphi-\varphi\right\vert \;\leq\;\left\vert V_{t}%
\varphi-S_{t}^{\alpha}\varphi\right\vert +\left\vert S_{t}^{\alpha}%
\varphi-\varphi\right\vert .
\end{equation}
By Corollary \ref{C.scSalpha}, it suffices to deal with the first term at the
right hand side. By equation (\ref{cumulant-int}), we have to look at%
\begin{equation}
\Big|\int_{0}^{t}\!\mathrm{d}s\;S_{t-s}^{\alpha}v^{\beta}(s)\,v(s)\Big|.
\end{equation}
But from domination (\ref{dom}) and Corollary \ref{Cor.flow.bound},%
\begin{equation}
0\;\leq\;v^{\beta}(s)\;\leq\;C\,(1+s^{-\kappa})\,\phi^{\beta},\qquad0<s\leq T,
\label{combined}%
\end{equation}
with $\,\kappa$\thinspace\ from (\ref{kappa}) and a constant $\,C=C(T)$%
\thinspace\ (note that other dependencies are not important in the present
proof). Thus, we can apply Lemma \ref{L.cl-0} to finish the proof.\hfill
\rule{0.5em}{0.5em}

\subsection{Uniqueness of solutions\label{subsection-uniqueness}}

The following lemma will be useful when we estimate the difference of
solutions to (\ref{cumulant-int}).

\begin{lemma}
[\textbf{An elementary observation}]\label{trivial} Let $\,\beta>0$%
\thinspace\ and $\,a,b\in\mathsf{R}$. Then
\begin{equation}
\left\vert _{\!_{\!_{\,}}}a\,(a\vee0)^{\beta}-b\,(b\vee0)^{\beta}\right\vert
\;\leq\;(1+\beta)\left(  _{\!_{\!_{\,}}}|a|+|b|\right)  ^{\beta}\,|a-b|.
\label{triv}%
\end{equation}

\end{lemma}

\proof First assume that $\,a,b\geq0.$\thinspace\ By the mean value theorem,
there exists a number $\,c$\thinspace\ between $\,a$\thinspace\ and
$\,b$\thinspace\ such that
\begin{equation}
\left|  a^{1+\beta}-b^{1+\beta}\right|  \;=\;(1+\beta)\,c^{\beta}%
\,|a-b|\;\leq\;(1+\beta)\,(a+b)^{\beta}\,|a-b|.
\end{equation}
This proves (\ref{triv}) for $a,b\geq0$.

Now suppose that $a,b<0$. In that case the left hand side in (\ref{triv})
disappears, hence (\ref{triv}) holds trivially.

Finally, it remains to consider the case $\,a<0\leq b$. Then,
\[
\left\vert a\,(a\vee0)^{\beta}-b\,(b\vee0)^{\beta}\right\vert \;=\;b^{1+\beta
}\;\leq\;(1+\beta)\,b^{\beta}\,b\;\leq\;(1+\beta)\left(  _{\!_{\!_{\,}}%
}|a|+|b|\right)  ^{\beta}\,|a-b|,
\]
and the proof is finished.\hfill\rule{0.5em}{0.5em}\bigskip

We are ready to prove uniqueness for solutions to (\ref{cumulant-int}).

\begin{lemma}
[Uniqueness]\label{unique-3}Fix $\,\varphi\in\Phi.$\thinspace\ Suppose that
$\,\,u,v$\thinspace\ are $\Phi$--valued solutions of equation\/
\emph{(\ref{cumulant-int})}. Then $\,u=v$.
\end{lemma}

%

\proof
Let
\begin{equation}
D(t,x)\;:=\;u(t,x)-v(t,x),\qquad t>0,\,\ x\neq0.
\end{equation}
Note that by Lemma \ref{eta-zero}, for $\,T>0$\thinspace\ fixed, $\,$%
\begin{equation}
\left\Vert D(t)_{\!_{\!_{\,}}}\right\Vert _{\mathcal{H}}\;\leq
\;2\,C_{\mathrm{\ref{eta-zero}}}\,\Vert\varphi\Vert_{\mathcal{H\,}}%
,\qquad0\leq t\leq T.\label{sup}%
\end{equation}
By the elementary inequality (\ref{triv}),
\begin{align}
\left\vert D(t,x)_{\!_{\!_{\,}}}\right\vert \; &  =\;\eta\left\vert \int
_{0}^{t}\mathrm{d}s\;S_{t-s}^{\alpha}\big(u^{1+\beta}(s)-v^{1+\beta
}(s)\big )(x)\right\vert \label{3-est-1}\\
\; &  \leq\;\eta\int_{0}^{t}\mathrm{d}s\;S_{t-s}^{\alpha}\big|u^{1+\beta
}(s)-v^{1+\beta}(s)\big|(x)\nonumber\\
\; &  \leq\;2\eta\int_{0}^{t}\mathrm{d}s\;S_{t-s}^{\alpha}{\Big(}%
\!\big [u^{\beta}(s)+v^{\beta}(s)\big]\left\vert D(s)_{\!_{\!_{\,}}%
}\right\vert \!{\Big)}(x).\nonumber
\end{align}
\textrm{F}rom (\ref{combined}), we get%
\begin{equation}
\left\vert D(t,x)_{\!_{\!_{\,}}}\right\vert \;\leq\;C\,\int_{0}^{t}%
\mathrm{d}s\;(1+s^{-\kappa})\,S_{t-s}^{\alpha}\left(  \left\vert
D(s)_{\!_{\!_{\,}}}\right\vert \phi^{\beta}\right)  (x).
\end{equation}
Thus%
\begin{gather}
\left\Vert D(t)_{\!_{\!_{\,}}}\right\Vert _{\mathcal{H}}\;\leq\;C\int_{0}%
^{t}\mathrm{d}s\;(1+s^{-\kappa})\,%
\Big\|%
S_{t-s}^{\alpha}\left(  \left\vert D(s)_{\!_{\!_{\,}}}\right\vert \phi^{\beta
}\right)
\Big\|%
_{\mathcal{H}}\nonumber\\
\leq\;C\int_{0}^{t}\mathrm{d}s\;(1+s^{-\kappa})\,C_{\mathrm{\ref{Cor.est.sing}%
}}\,(t-s)^{-\lambda}\left\Vert D(s)_{\!_{\!_{\,}}}\right\Vert _{\mathcal{H\,}%
},
\end{gather}
$0<t\leq T,$\thinspace\ where we used Corollary \ref{Cor.est.sing} and
notation (\ref{lambda}). Setting
\begin{equation}
D_{t}\;:=\;\sup_{0<s\leq t}\left\Vert D(s)_{\!_{\!_{\,}}}\right\Vert
_{\mathcal{H\,}},\qquad0<t\leq T,
\end{equation}
[for finiteness, recall (\ref{sup})], since $\,I(t)$\thinspace\ from
(\ref{not.I.t}) is increasing in $\,t$\thinspace\ [recall representation
(\ref{scal.I})], we find%
\begin{equation}
D_{t}\;\leq\;C\,D_{t}\,I(t),\qquad0<t\leq T,\label{found}%
\end{equation}
with some constant $\,C=C(T).$\thinspace\ Therefore, by (\ref{not.I.t}),
$\,D_{t}=0$\thinspace\ for all sufficiently small $\,t.$\thinspace\ Since the
model is time-homogeneous, we can repeat the argument finitely often to extent
to the whole interval $\,[0,T].$\thinspace\ (For a more complicated
time-inhomogeneous situation, see the proof of Lemma \ref{L.cl-1} below.)
Because $\,T$\thinspace\ is arbitrary, and since $\,u$\thinspace\ and
$\,v$\thinspace\ are $\,\Phi$--valued, we found $\,u=v,$\thinspace\ and the
proof is complete.\hfill\rule{0.5em}{0.5em}

\subsection{Auxiliary functions $v_{N,n}$}

For fixed integer $\,N\geq2$\thinspace\ set%
\begin{equation}
\psi_{N}\;:=\;\psi\wedge N. \label{not.psi.N}%
\end{equation}
For the fixed $\,T,N,\varphi,\psi,$\thinspace\ we inductively define functions
$\,v_{N,n\,}.$\thinspace\ First of all,%
\begin{equation}
v_{N,0}(t)\;:=\;S_{t}^{\alpha}\varphi\;\in\;\Phi,\qquad0\leq t\leq T.
\label{ind-1}%
\end{equation}
Assuming that we have defined $\,v_{N,n\,}$\thinspace\ for some $\,n,$%
\thinspace\ let%
\begin{equation}
v_{N,n+1}(t,x)\;:=\;S_{t}^{\alpha}\varphi\,(x)-\,\int_{0}^{t}\!\mathrm{d}%
s\;S_{t-s}^{\alpha}\!\left(  \psi_{N}(s)\,v_{N,n}(s)_{\!_{\!_{\,}}}\right)
(x), \label{ind-2}%
\end{equation}
$0\leq t\leq T,\;\,x\in\mathsf{\dot{R}}^{d},$\thinspace\ provided the latter
integral makes sense.

\begin{lemma}
[\textbf{Properties of }$v_{N,n}$]\label{L.cl-1}For all\/ $\,n\geq0$%
\thinspace\ and $\,t\in\lbrack0,T],$\thinspace\
\begin{equation}
0\;\leq\;v_{N,n}(t)\;\leq\;S_{t}^{\alpha}\varphi,\quad\text{and}\quad x\mapsto
v_{N,n}(t,x)\quad\text{is continuous.} \label{cl-1}%
\end{equation}

\end{lemma}

%

\proof
For $\,n=0,$\thinspace\ the claim is true by (\ref{ind-1}). Suppose that we
have verified (\ref{cl-1}) for some $\,n\geq0.$\thinspace\ Then the integral
in (\ref{ind-2}) is non-negative, hence%
\begin{equation}
v_{N,n+1}(t)\;\leq\;S_{t}^{\alpha}\varphi,\qquad t\in\lbrack0,T].
\end{equation}
Assume for the moment that $\,v_{N,n+1}\geq0$\thinspace\ under our induction
hypothesis. Then by Lemma \ref{L.cl-0}, $\,$%
\begin{equation}
v_{N,n+1}(t)\in\Phi,\qquad t\in\lbrack0,T],
\end{equation}
and the proof would be finished.

Next we will verify that $\,v_{N,n+1}$\thinspace\ is non-negative on
$\,[0,1/N].$\thinspace\ Since $\,\psi_{N}\leq N,$\thinspace\ and using the
induction assumption, it follows that%
\begin{equation}
S_{t-s}^{\alpha}\!\left(  \psi_{N}(s)\,v_{N,n}(s)_{\!_{\!_{\,}}}\right)
\;\leq\;S_{t-s}^{\alpha}\left(  N\,S_{s}^{\alpha}\varphi\right)
\;\leq\;N\,S_{t}^{\alpha}\varphi.
\end{equation}
Therefore, if $\,0\leq t\leq1/N,$%
\begin{equation}
v_{N,n+1}(t)\;\geq\;S_{t}^{\alpha}\varphi\,-\,N\int_{0}^{1/N}\!\mathrm{d}%
s\;S_{t}^{\alpha}\varphi\;=\;0.
\end{equation}

Now we prove that $\,v_{N,n+1}$\thinspace\ is non-negative on $\,[0,T].$%
\thinspace\ We use induction on the time intervals $\,\left[
k/N,\,(k+1)/N_{\!_{\!_{\,}}}\right]  \!,$\thinspace\ $0\leq k<NT.$%
\thinspace\ To begin with, we have already shown that $\,v_{N,n+1}$%
\thinspace\ is non-negative on $\,[0,1/N].$\thinspace\ Also, we know already
that $\,v_{N,n+1}(1/N)\in\Phi.$\thinspace\ Suppose that we have shown
$\,v_{N,n+1}$\thinspace\ is non-negative on $\,\left[  _{\!_{\!_{\,}}%
}(k-1)/N,\,k/N\right]  $\thinspace\ for some $\,0\leq k<NT-1,$\thinspace\ and
that $\,v_{N,n+1}(k/N)\in\Phi.$\thinspace\ We will shift time, and define%
\begin{subequations}
\begin{align}
v_{N,n+1}^{(k)}(t)\;  &  :=\;v_{N,n+1}(t+k/N),\\
\varphi_{N,n+1}^{(k)}(t)\;  &  :=\;v_{N,n+1}(k/N),\\
\psi_{N}^{(k)}(t)\;  &  :=\;\psi_{N}(t+k/N),
\end{align}
$0\leq t\leq1/N.$\thinspace\ We claim that%
\end{subequations}
\begin{equation}
v_{N,n+1}^{(k)}(t)\;:=\;S_{t}^{\alpha}\varphi_{N,n+1}^{(k)}-\int_{0}%
^{t}\!\mathrm{d}s\;S_{t-s}^{\alpha}\big(\psi_{N}^{(k)}(s)\,v_{N,n}%
^{(k)}(s)\big),\quad\;0\leq t\leq\frac{1}{N}\,. \label{shift-1}%
\end{equation}
Assume for the moment that (\ref{shift-1}) is true. Then the proof that
$\,v_{N,n+1}\geq0$\thinspace\ on $\,\left[  k/N,\,(k+1)/N_{\!_{\!_{\,}}%
}\right]  $\thinspace\ reduces to showing that $\,v_{N,n+1}^{(k)}\geq
0$\thinspace\ on $\,[0,1/N].$\thinspace\ But this follows from the step we
have already done. We are left with showing (\ref{shift-1}).

Using definition (\ref{ind-2}), we get%
\begin{equation}
S_{t}^{\alpha}v_{N,n+1}(r)\;:=\;S_{t+r}^{\alpha}\varphi\,-\,\int_{0}%
^{r}\!\mathrm{d}s\;S_{t+r-s}^{\alpha}\!\left(  \psi_{N}(s)\,v_{N,n}%
(s)_{\!_{\!_{\,}}}\right)  \!.\label{composition-1}%
\end{equation}
Let $\,r=k/N.$\ Then, for $\,0\leq t\leq1/N,$\
\begin{equation}
S_{t}^{\alpha}\varphi_{N,n+1}^{(k)}\;=\;S_{t+k/N}^{\alpha}\varphi\,-\,\int
_{0}^{k/N}\!\mathrm{d}s\;S_{t+k/N-s}^{\alpha}\!\left(  \psi_{N}(s)\,v_{N,n}%
(s)_{\!_{\!_{\,}}}\right)  \!.\label{composition-1a}%
\end{equation}
Also, by a change of variables, for $\,0\leq t\leq1/N,$\thinspace\ we obtain%
\begin{align}
&  \int_{0}^{t}\!\mathrm{d}s\;S_{t-s}^{\alpha}\!\left(  \psi_{N}%
^{(k)}(s)\,v_{N,n}^{(k)}(s)_{\!_{\!_{\,}}}\right)  \label{composition-2}%
\\[1pt]
&  =\;\int_{k/N}^{t+k/N}\!\mathrm{d}s\;S_{t+k/N-s}^{\alpha}\!\left(  \psi
_{N}(s)\,v_{N,n}(s)_{\!_{\!_{\,}}}\right)  \!.\nonumber
\end{align}
Inserting (\ref{composition-1a}) and (\ref{composition-2}) into the right hand
side of (\ref{shift-1}), and then using (\ref{ind-2}), we get%
\begin{align}
&  S_{t}^{\alpha}\varphi_{N,n+1}^{(k)}\,-\,\int_{0}^{t}\!\mathrm{d}%
s\;S_{t-s}^{\alpha}\!\big(\psi_{N}^{(k)}(s)\,v_{N,n}^{(k)}(s)\big)\\
&  =\;S_{t+k/N}^{\alpha}\varphi\,-\,\int_{0}^{t+k/N}\!\mathrm{d}%
s\;S_{t+k/N-s}^{\alpha}\!\left(  \psi_{N}(s)\,v_{N,n}^{{}}(s)_{\!_{\!_{\,}}%
}\right)  \\[2pt]
&  =\;v_{N,n+1}(t+k/N)\;=\;v_{N,n+1}^{(k)}(t),\nonumber
\end{align}
which proves (\ref{shift-1}). This finishes the proof.\hfill
\rule{0.5em}{0.5em}

\subsection{Auxiliary functions $v_{n}$}

Recall that we fixed $\,\varphi\in\Phi_{+\,}.$\thinspace\ For $\,n\geq
0,$\thinspace\ we inductively define functions $\,v_{n}$\thinspace\ as
follows. Let%
\begin{equation}
v_{0}(t)\;:=\;S_{t}^{\alpha}\varphi,\qquad0\leq t\leq T, \label{ind-3}%
\end{equation}
and, given $\,v_{n\,},$\thinspace\ set%
\begin{equation}
v_{n+1}(t,x)\;:=\;S_{t}^{\alpha}\varphi\,(x)\,-\,\int_{0}^{t}\!\mathrm{d}%
s\;S_{t-s}^{\alpha}\!\left(  \psi(s)\,v_{n}(s)_{\!_{\!_{\,}}}\right)  \!(x),
\label{ind-4}%
\end{equation}
$0\leq t\leq T,$\thinspace\ $x\in\mathsf{\dot{R}}^{d}.$

\begin{lemma}
[\textbf{Properties of }$v_{n}$]\label{L.cl-2}For each $\,n\geq0$%
\thinspace\ and $\,t\in\lbrack0,T],$%
\begin{equation}
\lim_{N\uparrow\infty}\,\left\vert v_{N,n}(t)-v_{n}(t)_{\!_{\!_{\,}}%
}\right\vert \;=\;0. \label{cl-2'}%
\end{equation}
Moreover,%
\begin{equation}
0\,\leq\,v_{n}(t)\,\leq\,S_{t}^{\alpha}\varphi,\quad\text{and}\quad x\mapsto
v_{n}(t,x)\quad\text{is continuous.} \label{cl-2}%
\end{equation}

\end{lemma}

%

\proof
Again, we use induction on $\,n.$\thinspace\ The claims are trivially true for
$\,n=0.$\thinspace\ Suppose they hold for $\,n.$\thinspace\ By definitions
(\ref{ind-2}) and (\ref{ind-4}),%
\begin{align}
&  \left\vert v_{n+1}(t)-v_{N,n+1}(t)_{\!_{\!_{\,}}}\right\vert
\;\label{cl-2-1}\\[1pt]
&  =\;\Big|\int_{0}^{t}\!\mathrm{d}s\;S_{t-s}^{\alpha}\!\left(  \psi
(s)\,v_{n}(s)-\psi_{N}(s)\,v_{N,n}(s)_{\!_{\!_{\,}}}\right)
\!\Big|\nonumber\\
&  \leq\;\bigg|\int_{0}^{t}\!\mathrm{d}s\;S_{t-s}^{\alpha}\!\left(  \left[
\psi(s)-\psi_{N}(s)_{\!_{\!_{\,}}}\right]  v_{n}(s)_{\!_{\!_{\,_{{}}}}%
}\right)  \!\bigg|\nonumber\\
&  \;\qquad+\;\int_{0}^{t}\!\mathrm{d}s\;S_{t-s}^{\alpha}\!\left(
\psi(s)\left\vert v_{n}(s)-v_{N,n}(s)_{\!_{\!_{\,}}}\right\vert
_{\!_{\!_{\,_{{}}}}}\right)  \ =:\;A_{N,n}+B_{N,n}\nonumber
\end{align}
with the obvious correspondence. We will show that both $\,A_{N+1,n}\,\ $and
$\,B_{N+1,n}$\thinspace\ tend to $\,0$\thinspace\ as $\,N\uparrow\infty
,$\thinspace\ giving (\ref{cl-2'}) for $\,n+1.$\thinspace\ This then yields
also the remaining claims in Lemma \ref{L.cl-2} for $\,n+1.$ In fact, by Lemma
\ref{L.cl-1} then the inequalities hold in (\ref{cl-2}), and Lemma
\ref{L.cl-0} gives the continuity claim.

First note that by the induction hypothesis,%
\begin{align}
A_{N,n}\;  &  \leq\;\int_{0}^{t}\!\mathrm{d}s\;S_{t-s}^{\alpha}\!\left(
\left[  \psi(s)-\psi_{N}(s)_{\!_{\!_{\,}}}\right]  S_{s}^{\alpha}%
\varphi_{\!_{\!_{\,_{{}}}}}\right) \label{A.est}\\[1pt]
&  \leq\;\int_{0}^{t}\!\mathrm{d}s\;S_{t-s}^{\alpha}\!\left(  \psi
(s)\,S_{s}^{\alpha}\varphi_{\!_{\!_{\,}}}\right)  \;<\;\infty,\nonumber
\end{align}
by Lemma \ref{L.cl-0}. Thus, by monotone convergence,%
\begin{equation}
\lim_{N\uparrow\infty}\,A_{N,n}\;=\;0. \label{A-zero}%
\end{equation}
By Lemma \ref{L.cl-1} and the induction hypothesis,%
\begin{equation}
\left\vert v_{N,n}(t)-v_{n}(t)_{\!_{\!_{\,}}}\right\vert \;\leq\;2\,S_{t}%
^{\alpha}\varphi.
\end{equation}
Moreover, by the induction assumption, (\ref{cl-2'}) holds. Then, by Lemma
\ref{L.cl-0}, the dominated convergence theorem implies that%
\begin{equation}
\lim_{N\uparrow\infty}\,B_{N,n}\;=\;0,
\end{equation}
finishing the proof.\hfill\rule{0.5em}{0.5em}

\subsection{A linearized equation}

Next we show that $\,v_{n}$\thinspace\ converges as $\,n\uparrow\infty
$\thinspace\ to a solution of a linearized equation.

\begin{lemma}
[\textbf{Linearized equation}]\label{from-trotter}Fix again $\,T\geq
1,$\thinspace\ $\varphi\in\Phi,$\thinspace\ and $\,\psi:\mathsf{R}_{+}%
\times\mathsf{R}^{d}\rightarrow\mathsf{R}_{+}$\thinspace\ satisfying\/
\emph{(\ref{a-1}). }Then, for $\,0\leq t\leq T,$\thinspace\ the (non-negative)
limit
\begin{equation}
v(t)\;:=\;\lim_{n\uparrow\infty}\,v_{n}(t)
\end{equation}
exists pointwise and is the unique $\Phi$--valued solution to
\begin{equation}
v(t)\;=\;S_{t}^{\alpha}\varphi\,-\,\int_{0}^{t}\mathrm{d}s\;S_{t-s}^{\alpha
}\!\left(  \psi(s)\,v(s)_{\!_{\!_{\,}}}\right)  \!,\qquad0\leq t\leq T,
\label{linearized}%
\end{equation}
implying $\,$%
\begin{equation}
0\,\leq\,v(t)\,\leq\,S_{t}^{\alpha}\varphi,\qquad0\leq t\leq T. \label{cl-4}%
\end{equation}
\thinspace\ 
\end{lemma}

%

\proof
This follows from the usual Picard iteration argument, in along the same lines
as in our uniqueness proof (Lemma \ref{unique-3}). Continuity again follows
from Lemma \ref{L.cl-0}.\hfill\rule{0.5em}{0.5em}

\subsection{Existence of solutions}

Our next goal is to use Lemma \ref{from-trotter} to prove existence of a
$\Phi$--valued solution for equation (\ref{cumulant-int}). Hypothesis
\ref{H.choice}$\,$ is still in force.

\begin{lemma}
[Existence]\label{exist-3}To each $\,\varphi\in\Phi,$\thinspace\ there exists
a $\Phi$--valued solution $\,v=V\!\varphi$\thinspace\ to the log-Laplace
equation\/ \emph{(\ref{cumulant-int})}.
\end{lemma}

%

\proof
We want to construct a sequence of $\Phi$--valued functions $\,v_{n}%
$\thinspace\ satisfying%
\begin{equation}
v_{n}(t)\;\leq\;S_{t}^{\alpha}\varphi.\label{approximation-bound}%
\end{equation}
In fact, if $\,n=0,$\thinspace\ set $\,v_{0}:=S^{\alpha}\varphi.$%
\thinspace\ Assume that we have already defined $\,v_{n}$\thinspace\ for some
$\,n\geq0.$\thinspace\ Note that by Corollary \ref{Cor.flow.bound},%
\begin{equation}
\left\vert v_{n}^{\beta}(t)\right\vert \;\leq\;M\,(1+t^{-\kappa})\,\phi
^{\beta}\label{136}%
\end{equation}
with $\,\kappa$\thinspace\ from (\ref{kappa}). Let $\,v_{n+1}$\thinspace\ be
the unique $\Phi$--valued solution to
\begin{equation}
v_{n+1}(t,x)\;=\;S_{t}^{\alpha}\varphi\,(x)-\int_{0}^{t}\mathrm{d}%
s\;S_{t-s}^{\alpha}\!\left(  v_{n}^{\beta}(s)\,v_{n+1}(s)\right)
(x)\label{cumulant-int'}%
\end{equation}
according to Lemma \ref{from-trotter}, implying (\ref{approximation-bound})
for $\,n+1.$\thinspace\ Altogether, by induction we defined $\Phi$--valued
functions $\,v_{n}$\thinspace\ satisfying (\ref{cumulant-int'}),
(\ref{approximation-bound}), and (\ref{136}).

For $\,n\geq0,$\thinspace\ let
\begin{equation}
D_{n}\;:=\;v_{n+1}-v_{n\,}.
\end{equation}
Then, as in the proof of the uniqueness Lemma \ref{unique-3}, using Lemma
\ref{trivial}, for fixed $\,T>0,$\thinspace\ we find
\begin{equation}
\left\vert D_{n+1}(t)_{\!_{\!_{\,}}}\right\vert \;\leq\;C\int_{0}%
^{t}\mathrm{d}s\;(1+s^{-\kappa})\,S_{t-s}^{\alpha}\!\left(  \left\vert
D_{n}(s)_{\!_{\!_{\,}}}\right\vert \phi^{\beta}\right)  ,
\end{equation}
$0\leq t\leq T,$\thinspace\ with a constant $\,C=C(T),$\thinspace\ and
\begin{gather}
\left\Vert D_{n+1}(t)_{\!_{\!_{\,}}}\right\Vert _{\mathcal{H}}\;\leq
\;C\int_{0}^{t}\mathrm{d}s\;(1+s^{-\kappa})%
\Big\|%
S_{t-s}^{\alpha}\!\left(  \left\vert D_{n}(s)_{\!_{\!_{\,}}}\right\vert
\phi^{\beta}\right)  \!%
\Big\|%
_{\mathcal{H}}\nonumber\\
\leq\;C\int_{0}^{t}\mathrm{d}s\;(1+s^{-\kappa})\,C_{\mathrm{\ref{Cor.est.sing}%
}}\,(t-s)^{-\lambda}\left\Vert D_{n}(s)_{\!_{\!_{\,}}}\right\Vert
_{\mathcal{H\,}}.
\end{gather}
Setting
\begin{equation}
D_{n,t}\;:=\;\sup_{0\leq s\leq t}\left\Vert D_{n}(s)_{\!_{\!_{\,}}}\right\Vert
_{\mathcal{H\,}},\qquad0\leq t\leq T,
\end{equation}
we found that%
\begin{equation}
D_{n+1,t}\;\leq\;\varepsilon_{t}\,D_{n,t\,},
\end{equation}
where the $\,\varepsilon_{t}\,\ $are independent of $\,n,$\thinspace\ and
$\,\varepsilon_{t}\rightarrow0$ as $t\downarrow0.$\thinspace\ Thus, if our
$\,T>0$\thinspace\ is small enough, then there exists a constant
$\,0<\gamma<1$\thinspace\ such that if $\,0\leq t\leq T,$\thinspace\ then
\begin{equation}
D_{n+1,t}\;\leq\;\gamma\,D_{n,t\,},
\end{equation}
and so
\begin{equation}
D_{n,t}\;\leq\;\gamma^{n}\,D_{0,t\,}.
\end{equation}
Therefore, we can define
\begin{equation}
v(t,x)\;:=\;\sum_{n=0}^{\infty}D_{n}(t,x)\;=\;\lim_{n\uparrow\infty}%
\,v_{n}(t,x),\qquad t\leq0,\,\ x\neq0\label{def.v}%
\end{equation}
where the limit is taken in $\,\mathcal{H}_{+\,}.$\thinspace\ 

\textrm{F}rom our construction, it follows that%
\begin{equation}
0\;\leq\;v(t)\;\leq\;S_{t}^{\alpha}\varphi\quad\text{and}\quad\left\vert
v^{\beta}(t)\right\vert \;\leq\;M\,(1+t^{-\kappa})\,\phi^{\beta}%
.\label{limit.v}%
\end{equation}
Now we want to show that $v$ satisfies equation (\ref{cumulant-int}) for
$0\leq t\leq T$. We start from definition (\ref{cumulant-int'}). First, by
(\ref{def.v}),
\begin{equation}
\lim_{n\rightarrow\infty}\left\Vert v_{n+1}(t)-v(t)_{\!_{\!_{\,}}}\right\Vert
_{\mathcal{H}}=0.\label{vn.to.v}%
\end{equation}
As for the integral terms, we first note that for $a,b,c\geq0$, by
(\ref{triv}) we have
\begin{gather}
\left\vert ab^{\beta}-c^{1+\beta}\right\vert \leq\left\vert a^{1+\beta
}-c^{1+\beta}\right\vert +\left\vert b^{1+\beta}-c^{1+\beta}\right\vert
\label{triv'}\\[2pt]
\leq(1+\beta)(a+c)^{\beta}|a-c|+(1+\beta)(b+c)^{\beta}\,|b-c|.\nonumber
\end{gather}
Therefore, using the second part of (\ref{limit.v}), we have
\begin{align}
&
\Big\|%
\int_{0}^{t}\!\mathrm{d}s\;S_{t-s}^{\alpha}\!\left(  v_{n+1}(s)\,v_{n}^{\beta
}(s)\right)  \,-\,\int_{0}^{t}\!\mathrm{d}s\;S_{t-s}^{\alpha}v^{1+\beta}(s)%
\Big\|%
_{\mathcal{H}}\\
&  \leq\;\int_{0}^{t}\!\mathrm{d}s\;\left\Vert S_{t-s}^{\alpha}\left\vert
v_{n+1}(s)\,v_{n}^{\beta}(s)-v^{1+\beta}(s)\right\vert _{\!_{\!_{\,_{{}}}}%
}\right\Vert _{\mathcal{H}}\nonumber\\
&  \leq\;C\int_{0}^{t}\!\mathrm{d}s\;(1+s^{-\kappa})\bigg\|S_{t-s}^{\alpha
}\Big(\phi^{\beta}\left\vert v_{n+1}(s)-v(s)_{\!_{\!_{\,}}}\right\vert
+\phi^{\beta}\left\vert v_{n}(s)-v(s)_{\!_{\!_{\,}}}\right\vert
\Big)\bigg \|_{\mathcal{H\,}}.\nonumber
\end{align}
By Corollary \ref{Cor.est.sing}, this chain of inequalities can be continued
with
\begin{align}
&  \leq\;C\int_{0}^{t}\!\mathrm{d}s\;(1+s^{-\kappa})\,(t-s)^{-\lambda
}\left\Vert v_{n+1}(s)-v(s)_{\!_{\!_{\,}}}\right\Vert _{\mathcal{H}%
}\label{last}\\
&  \qquad+C\int_{0}^{t}\!\mathrm{d}s\;(1+s^{-\kappa})\,(t-s)^{-\lambda
}\left\Vert v_{n}(s)-v(s)_{\!_{\!_{\,}}}\right\Vert _{\mathcal{H}}%
\;\underset{n\uparrow\infty}{\longrightarrow}\;0,\nonumber
\end{align}
by (\ref{vn.to.v}), domination, Corollary \ref{C.scSalpha}, and dominated
convergence. Thus, $v$ satisfies (\ref{cumulant-int}) in $\,\mathcal{H}%
_{+\,},$\thinspace\ for $t\leq T$\thinspace\ and for sufficiently small
$\,T.$\thinspace\ By induction on intervals, as in the proof of Lemma
\ref{L.cl-1}, we extend the solution $\,v$\thinspace\ from $[0,T]$ to all
times. By Lemma \ref{L.cl-0}, the constructed solution $\,v$\thinspace\ is
$\Phi$--valued.\hfill\rule{0.5em}{0.5em}\bigskip

\noindent\emph{Completion of the proof of Theorem}\/ \ref{T.well}. $\,$ With
Lemmas \ref{unique-3} and \ref{exist-3}, we already proved the uniqueness and
existence claims, respectively. The semigroup property of $\,V$\thinspace
\ follows from uniqueness of solutions, and the strong continuity of
$\,V$\thinspace\ from Lemma~\ref{eta-zero}.\hfill\rule{0.5em}{0.5em}

\section{Construction of $\,X$\label{S.Constr.X}}

Having available the well-posedness of the log-Laplace equation (Theorem
\ref{T.well}), we now construct the desired process $\,X$\thinspace\ via a
Trotter product approach to the related log-Laplace semigroup. For this
purpose, we introduce an approximating log-Laplace equation separating
continuous state branching with index $\,1+\beta,$\thinspace\ and mass flow
according to $\,S^{\alpha}$\thinspace\ on alternate small time intervals
(Lemma \ref{L.vn-eq}). Main work consists in showing that its solutions
converge to the solutions to the true log-Laplace equation (Proposition
\ref{P.vn-converge}). Then, in the final subsection, all pieces can easily put
together to get our main result, the existence Theorem \ref{T.ex}.

\subsection{Approximating equation}

Fix $\,n\geq1$\thinspace\ and $\,\varphi\in\Phi.$\thinspace\ We inductively
define measurable functions $\,v_{n}$\thinspace\ on $\,\mathsf{R}_{+}%
\times\mathsf{\dot{R}}^{d}.$\thinspace\ First of all,%
\begin{equation}
v_{n}(0)\;:=\;S_{1/n}^{\alpha}\varphi.\label{def-vn-1}%
\end{equation}
Assume for the moment $\,v_{n}(\frac{k}{n})$\thinspace\ is defined for some
$\,k\geq0.$\thinspace\ For $\,\frac{k}{n}\leq t<\frac{k+1}{n},$\thinspace\ set%
\begin{equation}
v_{n}(t,x)\;:=\;\frac{v_{n}(\frac{k}{n},x)}{\left[  1+\eta\,\beta
\,v_{n}^{\beta}(\frac{k}{n},x)\,(t-\frac{k}{n})\right]  ^{1/\beta}}\,,\qquad
x\neq0.\label{def-vn-expl}%
\end{equation}
Note that%
\begin{equation}
\frac{\partial}{\partial t}v_{n}(t,x)\;=\;-\eta\,v_{n}^{1+\beta}%
(t,x)\quad\text{on}\quad\Big(\tfrac{k}{n},\tfrac{k+1}{n}\Big)\times
\mathsf{\dot{R}}^{d},\label{def-vn-2}%
\end{equation}
that $\,v_{n}(\frac{k}{n}+\,,\,x)=v_{n}(\frac{k}{n},x),$\thinspace\ and that
also the limit $\,v_{n}(\frac{k+1}{n}-\,,\,x)$\thinspace\ exists. Note that
(\ref{def-vn-expl}) gives the log-Laplace function of continuous state
branching with index $\,1+\beta.$\thinspace\ Put%
\begin{equation}
v_{n}\Big(\tfrac{k+1}{n}\,,\,x\Big)\;:=\;S_{1/n}^{\alpha}v_{n}\Big(\tfrac
{k+1}{n}-\,,\,\cdot\,\Big)(x),\qquad x\neq0.\label{def-vn-3}%
\end{equation}

\begin{lemma}
[\textbf{Approximating log-Laplace equation}]\label{L.vn-eq}The function
$\,v_{n}\geq0$\thinspace\ we have just defined satisfies%
\begin{equation}
v_{n}(t,x)\;=\;S_{\left(  1+[tn]\right)  /n}^{\alpha}\varphi(x)\,-\,\eta
\int_{0}^{t}\!\mathrm{d}s\;S_{\left(  [tn]-[sn]\right)  /n}^{\alpha}\!\left(
v_{n}^{1+\beta}(s)\right)  (x) \label{vn-eq}%
\end{equation}
on $\,\mathsf{R}_{+}\times\mathsf{\dot{R}}^{d}.$
\end{lemma}

%

\proof
Differentiating equation (\ref{vn-eq}) to $\,t\neq\frac{k}{n},$\thinspace
\ $k\geq0,$\thinspace\ gives the true statement (\ref{def-vn-2}). On the other
hand, for $\,t=\frac{k}{n},$\thinspace\ $k\geq0,$\thinspace\
\begin{equation}
v_{n}\big(\tfrac{k}{n}\big)\;=\;S_{(1+k)/n}^{\alpha}\varphi\,-\,\eta\sum
_{i=1}^{k}\Big[S_{\left(  k-(i-1)\right)  /n}^{\alpha}\int_{(i-1)/n}%
^{i/n}\!\mathrm{d}s\;v_{n}^{1+\beta}(s)\Big].
\end{equation}
By (\ref{def-vn-2}) and the fundamental theorem of calculus, the right hand
side of the latter equation equals $\,v_{n}(\frac{k}{n}),$\thinspace
\ finishing the proof.\hfill\rule{0.5em}{0.5em}\bigskip

Set $\,$%
\begin{equation}
t_{n}:=[tn]/n.
\end{equation}
Since $\,v_{n}$\thinspace\ is non-negative, from equation (\ref{vn-eq}) we get
the \emph{domination}%
\begin{equation}
0\;\leq\;v_{n}(t)\;\leq\;S_{1/n+t_{n}}^{\alpha}\varphi,\qquad t\geq0,\quad
n\geq1, \label{vn-bound-eq}%
\end{equation}
implying by Corollary \ref{Cor.est.sing} with $\,\beta=0,$%
\begin{equation}
\left\Vert v_{n}(t)_{\!_{\!_{\,}}}\right\Vert _{\mathcal{H}}\;\leq
\;C_{\mathrm{(\ref{vn.H})}}\,\Vert\varphi\Vert,\qquad0\leq t\leq T,\quad
n\geq1, \label{vn.H}%
\end{equation}
for each $\,T>0,$\thinspace\ and where $\,C_{\mathrm{(\ref{vn.H})}%
}=C_{\mathrm{(\ref{vn.H})}}(T).$\thinspace\ In particular, $\,v_{n}$%
\thinspace\ is $\,\mathcal{H}_{+}$--valued. Our aim is to show that the
$\,v_{n}$\thinspace\ converge to the unique solution to (\ref{cumulant-int}).
For this purpose, we will need the following estimate.

\begin{lemma}
[\textbf{Pointwise bound}]\label{L.pwb}Impose Hypothesis\/
\emph{\ref{H.choice}. }To each $\,\varphi\in\mathcal{H}_{+}$\thinspace\ and
$\,T>0,$\thinspace\ there is a $\,\varphi_{0}=\varphi_{0}(d,T,\alpha
,\varphi,\varrho)\,\ $in $\,\mathcal{H}_{+\,},$\thinspace\ such that%
\begin{equation}
\sup_{T/2\,\leq\,t\,\leq\,T}S_{t}^{\alpha}\varphi\;\leq\;\varphi
_{0\,}.\label{pwd}%
\end{equation}

\end{lemma}

%

\proof
Recall from (\ref{p-diff}) that%
\begin{equation}
P^{\alpha}(t;x,y)\;\leq\;P(t;x,y)+C_{\mathrm{\ref{L.kern.b}}}\,\bar{P}(t;x,y).
\end{equation}
Choose a constant $\,C_{\mathrm{(\ref{p.est})}}=C_{\mathrm{(\ref{p.est})}%
}(d,T)$\thinspace\ such that, for $\,T/2\leq t\leq T,$%
\begin{equation}
P(t;x,y)\;\leq\;C_{\mathrm{(\ref{p.est})}}P(T;x,y) \label{p.est}%
\end{equation}
and%
\begin{equation}
\bar{P}(t;x,y)\;=\;t^{-1/2}\,\phi(x)\,\phi(y)\,\mathrm{e}^{-|x|^{2}%
/4t}\,\,\mathrm{e}^{-|y|^{2}/4t}\;\leq\;C_{\mathrm{(\ref{p.est})}}\,\bar
{P}(T;x,y) \label{pbar.est}%
\end{equation}
[recall (\ref{bardef})]. \textrm{F}rom Lemma \ref{new-mink} (with $\,\beta
=0)$\thinspace\ we conclude that $\,S_{T}\varphi$\thinspace\ belongs to
$\,\mathcal{H}_{+\,},$\thinspace\ whereas Lemma \ref{new-mink.II} (with
$\,\beta=0)$\thinspace\ gives $\,\bar{S}_{T}\varphi\in\mathcal{H}_{+\,}%
.$\thinspace\ Therefore, we may set%
\begin{equation}
\varphi_{0}\;:=\;C_{\mathrm{(\ref{p.est})}}\left(  S_{T}\varphi
+C_{\mathrm{\ref{L.kern.b}}}\,\bar{S}_{T}\varphi\right)
\end{equation}
to finish the proof.\hfill\rule{0.5em}{0.5em}

\subsection{Convergence to the limit equation}

With the functions $\,v_{n}$\thinspace\ we may pass to the limit.

\begin{prop}
[\textbf{Convergence to the limit equation}]\label{P.vn-converge}Let
$\,\varphi\in\Phi.$\ Define $\,v_{n}$\thinspace\ as in\/ \emph{(\ref{def-vn-1}%
)--(\ref{def-vn-3}).} Let $\,v$\thinspace\ be the unique $\Phi$--valued
solution to\/ \emph{(\ref{cumulant-int})} according to Theorem\/
\emph{\ref{T.well}. }Then, for each $\,t\geq0,$
\begin{equation}
\lim_{n\uparrow\infty}\left\Vert v(t)-v_{n}(t)_{\!_{\!_{\,}}}\right\Vert
_{\mathcal{H}}\;=\;0.
\end{equation}

\end{prop}

%

\proof
We may restrict our attention to $\,t\in\lbrack0,T]$\thinspace\ for any
$\,T>1.$\thinspace\ \textrm{F}rom equations (\ref{cumulant-int}) and
(\ref{vn-eq}) we have%
\begin{align}
&  \left\Vert v(t)-v_{n}(t)_{\!_{\!_{\,}}}\right\Vert _{\mathcal{H}}%
\;\leq\;\left\Vert S_{t}^{\alpha}\varphi-S_{1/n+t_{n}}^{\alpha}\varphi
\right\Vert _{\mathcal{H}}\label{gronwall-vn}\\[2pt]
&  \qquad+\;\eta\int_{0}^{t_{n}}\!\mathrm{d}s\,\left\Vert S_{t-s}^{\alpha
}v^{1+\beta}(s)-S_{t_{n}-s_{n}}^{\alpha}v^{1+\beta}(s)\right\Vert
_{\mathcal{H}}\nonumber\\
&  \qquad+\;\eta\int_{0}^{t_{n}}\!\mathrm{d}s\;%
\Big\|%
S_{t_{n}-s_{n}}^{\alpha}\!\left\vert v^{1+\beta}(s)-v_{n}^{1+\beta
}(s)\right\vert \!%
\Big\|%
_{\mathcal{H}}\nonumber\\
&  \qquad+\;\eta\,\,%
\Big\|%
\int_{t_{n}}^{t}\!\mathrm{d}s\;S_{t-s}^{\alpha}v^{1+\beta}(s)%
\Big\|%
_{\mathcal{H}}+\eta\,\,%
\Big\|%
\int_{t_{n}}^{t}\!\mathrm{d}s\;S_{t_{n}-s_{n}}^{\alpha}v_{n}^{1+\beta}(s)%
\Big\|%
_{\mathcal{H}}\nonumber\\[2pt]
&  =:\;A_{n}(t)+B_{n}(t)+C_{n}(t)+D_{n}(t)+E_{n}(t)\nonumber
\end{align}
with the obvious correspondence. We will deal with each of these terms
separately.\medskip

\noindent1$^{\circ}$ $\left(  A_{n}(t)_{\!_{\!_{\,}}}\right)  \!.$%
\quad\textrm{F}rom the semigroup property and boundedness (\ref{bounded}),%
\begin{equation}
A_{n}(t)\;=\;\left\Vert S_{t}^{\alpha}\varphi-S_{1/n+t_{n}}^{\alpha}%
\varphi\right\Vert _{\mathcal{H}}\;\leq\;C_{\mathrm{\ref{Cor.est.sing}}%
}\,\left\Vert S_{\left\vert t-1/n-t_{n}\right\vert }^{\alpha}\varphi
-\varphi\right\Vert _{\mathcal{H}\,}.
\end{equation}
But%
\begin{equation}
\Big|t-\frac{1}{n}-t_{n}\Big|\,\leq\,\frac{2}{n}\,,\qquad t\geq0.\label{tn}%
\end{equation}
Hence,%
\begin{equation}
\sup_{0\leq t\leq T}A_{n}(t)\;\leq\;C_{\mathrm{\ref{Cor.est.sing}}}%
\,\sup_{0\leq s\leq2/n}\left\Vert S_{s}^{\alpha}\varphi-\varphi\right\Vert
_{\mathcal{H}\,}\;\underset{n\uparrow\infty}{\longrightarrow}\;0\label{term-a}%
\end{equation}
by strong continuity according to Corollary \ref{C.scSalpha}.\medskip

\noindent2$^{\circ}$ $\left(  D_{n}(t)_{\!_{\!_{\,}}}\right)  \!.$\quad
Clearly, by our estimates,%
\begin{equation}
D_{n}(t)\;\leq\;\eta\int_{t_{n}}^{t}\!\mathrm{d}s\;(1+s^{-\kappa
})\,(t-s)^{-\lambda}\,\Vert\varphi\Vert_{\mathcal{H}\,}.
\end{equation}
By scaling, the integral equals%
\[
t^{1-\lambda}\int_{t_{n}/t}^{1}\!\mathrm{d}s\;(1-s)^{-\lambda}+t^{1-\kappa
-\lambda}\int_{t_{n}/t}^{1}\!\mathrm{d}s\;s^{-\kappa}\,(1-s)^{-\lambda
}\;=:\;I_{n}(t)+I\!I_{n}(t)
\]
with the obvious correspondence. Take $\,\varepsilon\in(0,T)$\thinspace\ and
let $\,n>1/\varepsilon.$\thinspace\ Since%
\begin{equation}
\frac{t_{n}}{t}\;\geq\;1-\frac{1}{tn}\;\geq\;1-\frac{1}{\varepsilon
n}\,,\qquad\varepsilon\leq t\leq T,
\end{equation}
we get%
\begin{equation}
\sup_{\varepsilon\leq t\leq T}I_{n}(t)\;\leq\;T^{1-\lambda}\int
_{1-1/\varepsilon n}^{1}\!\mathrm{d}s\;(1-s)^{-\lambda}\;\underset
{n\uparrow\infty}{\longrightarrow}\;0,
\end{equation}
whereas%
\begin{equation}
\sup_{0\leq t\leq\varepsilon}I_{n}(t)\;\leq\;\varepsilon^{1-\lambda}\int
_{0}^{1}\!\mathrm{d}s\;(1-s)^{-\lambda}\;\underset{n\uparrow\infty
}{\longrightarrow}\;0.
\end{equation}
Consequently,%
\begin{equation}
\sup_{0\leq t\leq T}I_{n}(t)\;\underset{n\uparrow\infty}{\longrightarrow}\;0,
\end{equation}
and the same reasoning leads to the analogous statement on $\,I\!I_{n}%
(t).$\thinspace\ Summarizing,%
\begin{equation}
\sup_{0\leq t\leq T}D_{n}(t)\;\underset{n\uparrow\infty}{\longrightarrow
}\;0.\label{(D)}%
\end{equation}
\smallskip

\noindent3$^{\circ}$ $\left(  E_{n}(t)_{\!_{\!_{\,}}}\right)  \!.$\quad Assume
that $\,t>t_{n\,}.$\thinspace\ By (\ref{vn-eq}),
\begin{equation}
E_{n}(t)\;=\;\left\Vert v_{n}(t)-v_{n}(t_{n})_{\!_{\!_{\,}}}\right\Vert
_{\mathcal{H\,}}.
\end{equation}
According to the definition (\ref{def-vn-expl}) of $\,v_{n}(t),$%
\begin{align}
0\; &  \leq\;v_{n}(t_{n\,},x)-v_{n}(t,x)\;\nonumber\\[1pt]
&  =\;v_{n}(t_{n},x)\,\bigg(1-\frac{1}{\left[  1+\eta\,\beta\,v_{n}^{\beta
}(t_{n},x)\,(t-t_{n})\right]  ^{1/\beta}}\bigg).\label{But}%
\end{align}
Using domination and Lemma \ref{L.pwb}, there is a $\,\varphi_{0}=\varphi
_{0}(\varphi,t,\alpha,d,\varrho)\in\mathcal{H}_{+}$\thinspace\ such that%
\begin{equation}
v_{n}(t_{n\,},x)\;\leq\;S_{1/n+t_{n}}^{\alpha}\varphi\,(x)\;\leq
\;\varphi_{0\,},\label{v.phi}%
\end{equation}
since $\,t\leq1/n+t_{n}\leq1+t.$\thinspace\ But (\ref{But}) is increasing in
$\,v_{n}(t_{n\,},x),$\thinspace\ so we may insert (\ref{v.phi}) to obtain%
\[
0\;\leq\;v_{n}(t_{n\,},x)-v_{n}(t,x)\;\leq\;\varphi_{0}(x)\bigg(1-\frac
{1}{\big[1+\varphi_{0}^{\beta}(x)\big
]^{1/\beta}}\bigg)\;\leq\;\varphi_{0}(x),
\]
since $\,0\leq t-t_{n}\leq1/n.$\thinspace\ Then, from dominated convergence we
get%
\begin{equation}
\sup_{0\leq t\leq T}E_{n}(t)\;\underset{n\uparrow\infty}{\longrightarrow
}\;0.\label{E}%
\end{equation}
\smallskip

\noindent4$^{\circ}$ $\left(  B_{n}(t)_{\!_{\!_{\,}}}\right)  .$\quad First of
all, we want to deal with $\,B_{n}(t)$ for small $\,t.$\thinspace\ Clearly,%
\begin{equation}
\left\Vert S_{t-s}^{\alpha}v^{1+\beta}(s)\right\Vert _{\mathcal{H}}%
\;\leq\;C\,(1+s^{-\kappa})\,(t-s)^{-\lambda}\,\Vert\varphi\Vert_{\mathcal{H\,}%
}\label{clearly1}%
\end{equation}
and%
\begin{equation}
\left\Vert S_{t_{n}-s_{n}}^{\alpha}v^{1+\beta}(s)\right\Vert _{\mathcal{H}%
}\;\leq\;C\,(1+s^{-\kappa})\,(t_{n}-s)^{-\lambda}\,\Vert\varphi\Vert
_{\mathcal{H\,}}.\label{clearly2}%
\end{equation}
Let $\,0<\varepsilon<T.$\thinspace\ Using notation (\ref{not.I.t}), from
(\ref{clearly1}) and (\ref{clearly2}), since $\,$%
\begin{equation}
t_{n}-s_{n}\,\geq\,t_{n}-s\,\geq\,0,\label{ts}%
\end{equation}
for $\,0\leq t\leq\varepsilon,$%
\begin{equation}
B_{n}(t)\;\leq\;C\left[  I(t)+I(t_{n})_{\!_{\!_{\,}}}\right]  \!.
\end{equation}
Moreover, since $\,I$\thinspace\ is increasing [recall (\ref{scal.I})],%
\begin{equation}
\sup_{n\geq1}\,\sup_{0\leq t\leq\varepsilon}B_{n}(t)\;\leq\;C\,\sup_{0\leq
t\leq\varepsilon}I(t)\;\underset{\varepsilon\downarrow0}{\longrightarrow
}\;0.\label{small.t}%
\end{equation}

Now we may restrict to $\,t\in\lbrack\varepsilon,T].$\thinspace\ We want to
exploit the strong continuity of the semigroup $\,S^{\alpha}$\thinspace
\ acting on $\,\mathcal{H}_{+}$\thinspace\ (Corollary \ref{C.scSalpha}). To
this end, we truncate $\,v^{1+\beta}$\thinspace\ to a function in
$\,\mathcal{H}_{+\,},$\thinspace\ and consider a small time interval around
$\,t_{n}$\thinspace\ separately to get rid of the varying upper integration
bound. Here are the details.

Take $\,\delta\in(0,\varepsilon)$\thinspace\ and $\,N\geq1.$\thinspace\ Set%
\begin{subequations}
\begin{align}
v_{1,N}(t)\;  &  :=\;\left(  v(t)\wedge N_{\!_{\!_{\,}}}\right)
\mathsf{1}_{B_{N}(0)}\\[1pt]
v_{2,N}(t)\;  &  :=\;v(t)-v_{1,N}(t).
\end{align}
Then for $\,\varepsilon\leq t\leq T,$%
\end{subequations}
\begin{align}
B_{n}(t)\;  &  \leq\;\int_{0}^{t-\delta}\!\mathrm{d}s\,\left\|  S_{t-s}%
^{\alpha}v_{1,N}^{1+\beta}(s)-S_{t_{n}-s_{n}}^{\alpha}v_{1,N}^{1+\beta
}(s)\right\|  _{\mathcal{H}}\\
&  \;+\ \int_{t-\delta}^{t}\!\mathrm{d}s\,\left\|  S_{t-s}^{\alpha}v^{1+\beta
}(s)\right\|  _{\mathcal{H}}+\int_{t-\delta}^{t_{n}}\!\mathrm{d}s\,\left\|
S_{t_{n}-s_{n}}^{\alpha}v^{1+\beta}(s)\right\|  _{\mathcal{H}}\nonumber\\
&  \;+\ \int_{0}^{t-\delta}\!\mathrm{d}s\,\left\|  S_{t-s}^{\alpha}%
v_{2,N}^{1+\beta}(s)\right\|  _{\mathcal{H}}+\int_{0}^{t-\delta}%
\!\mathrm{d}s\,\left\|  S_{t_{n}-s_{n}}^{\alpha}v_{2,N}^{1+\beta}(s)\right\|
_{\mathcal{H}}\nonumber\\[2pt]
&  =:\;B_{n}^{(1)}(t)+\cdots+B_{n}^{(5)}(t)\nonumber
\end{align}
in the obvious correspondence. Again, we deal with all terms
separately.\medskip

\noindent4.1$^{\circ}$ $\big(B_{n}^{(2)}(t)\big).$\quad\textrm{F}rom
(\ref{clearly1}) and scaling,%
\begin{equation}
B_{n}^{(2)}(t)\;\leq\;C\int_{1-\delta/t}^{1}\!\mathrm{d}s\;(1-s)^{-\lambda
}+C\int_{1-\delta/t}^{1}\!\mathrm{d}s\;s^{-\kappa}\,(1-s)^{-\lambda}.
\end{equation}
But $\,1-\delta/t\geq1-\delta/\varepsilon,$\thinspace\ hence%
\begin{equation}
\sup_{n}\,\sup_{\varepsilon\leq t\leq T}B_{n}^{(2)}(t)\;\underset
{\delta\downarrow0}{\longrightarrow}\;0.\label{B2}%
\end{equation}
\smallskip

\noindent4.2$^{\circ}$ $\big(B_{n}^{(3)}(t)\big).$\quad Similarly, from
(\ref{clearly2}) and (\ref{ts}),%
\begin{equation}
B_{n}^{(3)}(t)\;\leq\;C\int_{t-\delta}^{t_{n}}\!\mathrm{d}s\;(1+s^{-\kappa
})\,(t_{n}-s)^{-\lambda}.
\end{equation}
Now $\,$%
\begin{equation}
t\geq t_{n}\geq t-\frac{1}{n}\geq\varepsilon-\frac{1}{n}\geq\delta\label{now}%
\end{equation}
provided that $\,n\geq1/(\varepsilon-\delta).$\thinspace\ Thus, the lower
integration bound can be replaced by $\,t_{n}-\delta,$\thinspace\ and by
scaling,%
\begin{equation}
B_{n}^{(3)}(t)\;\leq\;C\int_{1-\delta/t_{n}}^{1}\!\mathrm{d}s\;(1-s)^{-\lambda
}+C\int_{1-\delta/t_{n}}^{1}\!\mathrm{d}s\;s^{-\kappa}\,(1-s)^{-\lambda}.
\end{equation}
By (\ref{now}), the lower integration bounds can be changed to $\,1-\delta
/(\varepsilon-1/n),$\thinspace\ implying%
\begin{equation}
\limsup_{n\uparrow\infty}\,\sup_{\varepsilon\leq t\leq T}B_{n}^{(3)}%
(t)\;\underset{\delta\downarrow0}{\longrightarrow}\;0.\label{B3}%
\end{equation}
\smallskip

\noindent4.3$^{\circ}$ $\big(B_{n}^{(4)}(t)\;$and$\;B_{n}^{(5)}(t)\big).$\quad
First note that for $\,s\in\lbrack0,t-\delta],$%
\begin{equation}
t-s\quad\text{and}\quad t_{n}-s_{n}\quad\text{belong to}\quad\lbrack
\delta/2,t\}\quad\text{if}\quad n>2/\delta\label{both}%
\end{equation}
(for instance, $\,t_{n}-s_{n}\geq t-\frac{1}{n}-s\geq-\frac{1}{n}+\delta
\geq\delta\,).$\thinspace\ Then, by domination and Corollary
\ref{Cor.flow.bound},%
\begin{equation}
v_{2,N}^{\beta}(s)\;\leq\;v^{\beta}(s)\;\leq\;(S_{s}^{\alpha}\varphi)^{\beta
}\;\leq\;C\,(1+s)^{-\kappa}\,\phi^{\beta},
\end{equation}
hence, by Lemma \ref{new-mink.II}, for $\,r\in\lbrack\delta/2,t]$%
\thinspace\ and $\,n>2/\delta,s\leq T,$%
\[
\big\|S_{r}^{\alpha}v_{2,N}^{1+\beta}(s)\big\|_{\mathcal{H}}\;\leq
\;C\,(1+s)^{-\kappa}\left\Vert S_{r}^{\alpha}\big(v_{2,N}(s)\,\phi^{\beta
}\big
)\right\Vert _{\mathcal{H}}\;\leq\;C\,(1+s)^{-\kappa}\left\Vert v_{2,N}%
(s)\right\Vert _{\mathcal{H}\,}.
\]
Now by domination and boundedness, for $\,0\leq s\leq T,$%
\begin{equation}
\left\Vert v_{2,N}(s)\right\Vert _{\mathcal{H}}\;\leq\;\left\Vert
S_{s}^{\alpha}\varphi\right\Vert _{\mathcal{H}}\;\leq\;\left\Vert
\varphi\right\Vert _{\mathcal{H}\,}%
\end{equation}
and%
\begin{equation}
\left\Vert v_{2,N}(s)\right\Vert _{\mathcal{H}}\;\underset{N\uparrow\infty
}{\searrow}\;0.
\end{equation}
Therefore,%
\begin{align}
&  \limsup_{n\uparrow\infty}\sup_{\varepsilon\leq t\leq T}\big(B_{n}%
^{(4)}(t)+B_{n}^{(5)}(t)\big)\label{B45}\\
&  \leq\;C\!\int_{0}^{T}\!\mathrm{d}s\;(1+s)^{-\kappa}\left\Vert
v_{2,N}(s)\right\Vert _{\mathcal{H}}\,\underset{N\uparrow\infty}{\searrow
}\,0,\nonumber
\end{align}
by monotone convergence, for the fixed $\,\varepsilon$\thinspace\ and
$\,\delta.$\medskip

\noindent4.4$^{\circ}$ $\big(B_{n}^{(1)}(t)\big).$\quad It remains to deal
with $\,B_{n}^{(1)}(t).$\thinspace\ By the semigroup property, boundedness,
and strong continuity,%
\begin{align}
&  \left\Vert S_{t-s}^{\alpha}v_{1,N}^{1+\beta}(s)-S_{t_{n}-s_{n}}^{\alpha
}v_{1,N}^{1+\beta}(s)\right\Vert _{\mathcal{H}}\nonumber\\[2pt]
&  \leq\;\sup_{0\leq r\leq2/n}\,\left\Vert S_{r}^{\alpha}v_{1,N}^{1+\beta
}(s)-v_{1,N}^{1+\beta}(s)\right\Vert _{\mathcal{H}}\;\underset{N\uparrow
\infty}{\searrow}\;0\label{suppp}%
\end{align}
for all $\,s$\thinspace\ and $\,N,$\thinspace\ since by definition%
\begin{equation}
v_{1,N}^{1+\beta}(s)\;\leq\;N^{1+\beta}\,\mathsf{1}_{B_{N}(0)}\;\in
\;\mathcal{H}_{+\,}.
\end{equation}
Moreover, the supremum in (\ref{suppp}) is bounded form above by%
\begin{equation}
2\,N^{1+\beta}\,\big\|\mathsf{1}_{B_{N}(0)}\big\|_{\mathcal{H}\,}.
\end{equation}
Therefore,%
\begin{equation}
\sup_{\varepsilon\leq t\leq T}B_{n}^{(1)}(t)\;\leq\;\int_{0}^{T}%
\!\mathrm{d}s\,\sup_{0\leq r\leq2/n}\,\left\Vert S_{r}^{\alpha}v_{1,N}%
^{1+\beta}(s)-v_{1,N}^{1+\beta}(s)\right\Vert _{\mathcal{H}}\;\underset
{n\uparrow\infty}{\searrow}\;0,\label{B1}%
\end{equation}
by monotone convergence, for all our $\,N,\varepsilon,\delta.$\medskip

\noindent4.5$^{\circ}$ (\emph{Conclusion})$.$\quad Putting together
(\ref{small.t}), (\ref{B2}), (\ref{B3}), (\ref{B45}), and (\ref{B1}),
\begin{equation}
\sup_{0\leq t\leq T}B_{n}(t)\;\underset{n\uparrow\infty}{\longrightarrow
}\;0.\label{B}%
\end{equation}

\noindent5$^{\circ}$ $\left(  C_{n}(t)_{\!_{\!_{\,}}}\right)  \!.$\quad First
note that%
\begin{equation}
C_{n}(t)\,=\,0\quad\text{for}\quad t\leq1/n. \label{0}%
\end{equation}
So we may assume that $\,t\geq1/n.$\thinspace\ Next we apply Lemma
\ref{trivial} to get for the term in abstract value sign in the definition of
$\,C_{n}(t)$\thinspace\ the bound%
\begin{equation}
C\left[  |v|^{\beta}(s)+|v_{n}|^{\beta}(s)_{\!_{\!_{\,}}}\right]  \left|
v(s)-v_{n}(s)_{\!_{\!_{\,}}}\right|  \!.
\end{equation}
\textrm{F}rom domination, the expression in square brackets is bounded by%
\begin{equation}
\left(  S_{s}^{\alpha}\varphi\right)  ^{\beta}+(S_{1/n+s_{n}}^{\alpha}%
\varphi)^{\beta}\;\leq\;C\,(1+s^{-\kappa})\,\phi^{\beta},
\end{equation}
where we used Corollary \ref{Cor.flow.bound} and $\,1/n+s_{n}\geq
s.$\thinspace\ But by Corollary \ref{Cor.est.sing},
\begin{equation}%
\Big\|%
S_{t_{n}-s_{n}}^{\alpha}\left|  v(s)-v_{n}(s)_{\!_{\!_{\,}}}\right|
\phi^{\beta}%
\Big\|%
_{\mathcal{H}}\;\leq\;C\,(t_{n}-s_{n})^{-\lambda}\left\|  v(s)-v_{n}%
(s)_{\!_{\!_{\,}}}\right\|  _{\mathcal{H}\,}.
\end{equation}
Setting%
\begin{equation}
F_{n}(t)\;:=\;\sup_{s\leq t}\,\left\|  v(s)-v_{n}(s)_{\!_{\!_{\,}}}\right\|
_{\mathcal{H}\,},
\end{equation}
we found for $\,\frac{1}{n}\leq t\leq T,$%
\begin{align}
&  \int_{0}^{t_{n}}\!\mathrm{d}s\,\,%
\Big\|%
S_{t_{n}-s_{n}}^{\alpha}\left|  v^{1+\beta}(s)-v_{n}^{1+\beta}(s)\right|
\Big\|%
_{\mathcal{H}}\nonumber\\
&  \leq\;C\,F_{n}(t)\int_{0}^{t_{n}}\!\mathrm{d}s\;(1+s^{-\kappa}%
)\,(t_{n}-s)^{-\lambda}\;\leq\;C\,F_{n}(t)\,I(t) \label{C2}%
\end{align}
[recall (\ref{not.I.t})].\medskip

\noindent6$^{\circ}$ (\emph{Completion of the proof}).\quad Inserting
(\ref{term-a}), (\ref{(D)}), (\ref{E}), (\ref{B}), (\ref{C2}) into
(\ref{gronwall-vn}), we obtain%
\begin{equation}
\left\Vert v(t)-v_{n}(t)_{\!_{\!_{\,}}}\right\Vert _{\mathcal{H}\,}%
\;\leq\;\varepsilon_{n}+C\,F_{n}(t)\,I(t),\label{obtain}%
\end{equation}
where $\,C=C(T)$\thinspace\ and where $\,\varepsilon_{n}=\varepsilon_{n}%
(T)$\thinspace\ is independent of $\,t$\thinspace\ and tends to $\,0$%
\thinspace\ as $\,n\uparrow\infty.$\thinspace\ Restrict further to
$\,t$\thinspace\ such that $\,C\,I(t)<\frac{1}{2}$\thinspace\ [recall
(\ref{not.I.t})] to get%
\begin{equation}
F_{n}(t)\;\leq\;\varepsilon_{n}+\frac{1}{2}\,F_{n}(t),\quad\text{that is}\quad
F_{n}(t)\;\leq\;2\,\varepsilon_{n}\;\underset{n\uparrow\infty}{\longrightarrow
}\;0.
\end{equation}
Consequently,%
\begin{equation}
\left\Vert v(t)-v_{n}(t)_{\!_{\!_{\,}}}\right\Vert _{\mathcal{H}\,}%
\;\underset{n\uparrow\infty}{\longrightarrow}\;0\quad\text{for all
sufficiently small}\,\ t.
\end{equation}
Repeating the argument, we can lift up for $\,t\in\lbrack0,T].$\thinspace
\ Since $\,T$\thinspace\ was arbitrary, the proof of Lemma \ref{P.vn-converge}
is finished altogether.\hfill\rule{0.5em}{0.5em}

\subsection{Construction of the process\label{SS.constrX}}

Here is now the more precise formulation of our main result, announced in
Theorem \ref{T.ex.0}:

\begin{theorem}
[\textbf{Existence of }$\,X$]\label{T.ex}Under Hypothesis\/
\emph{\ref{H.choice}, }there is a (unique in law) non-degen\-erate
$\mathcal{M}(\mathsf{\dot{R}}^{d})$--valued (time-homogeneous) Markov process
$\,X=(X,\,P_{\mu\,},\,\mu\in\mathcal{M})$\ with log-Laplace transition
functional\/ \emph{(\ref{logLap})}, where $\,v=V\!\varphi$\thinspace\ solves
equation\/ \emph{(\ref{cumulant-int})} with $\,\varphi\in\Phi$\thinspace
\ satisfying\/ \emph{(\ref{phi.bound}).}
\end{theorem}

\begin{remark}
[\textbf{Non-degeneracy}]\label{R.nond}It is easy to see that the following
expectation formula holds:%
\begin{equation}
P_{\mu}\!\left\langle X_{t\,},\varphi\right\rangle \ =\ \left\langle \mu
,S_{t}^{\alpha}\varphi\right\rangle ,\qquad\mu\in\mathcal{M},\quad
t\geq0,\quad\varphi\in\Phi.
\end{equation}
But $\,$%
\begin{equation}
V_{t}\varphi\ \neq\ S_{t}^{\alpha}\varphi,\qquad t>0,\quad\varphi\in\Phi
,\quad\varphi\neq0.
\end{equation}
Hence, the log-Laplace formula (\ref{logLap}) shows that $\,X$\thinspace\ is
different from its expectation, that is, it is non-degenerate.\hfill$\Diamond$
\end{remark}

\noindent\emph{Proof of Theorem}\/ \ref{T.ex}. $\,$Under Hypothesis
\ref{H.choice}, for $\,\mu\in\mathcal{M}=\mathcal{M}(\mathsf{\dot{R}}^{d}%
),$\thinspace\ from definitions (\ref{def-vn-1})--(\ref{def-vn-3}) of
$\,v_{n}=:V^{n}\varphi,$\thinspace\ $\varphi\in\Phi,$\thinspace\ for fixed
$\,t>0$\thinspace\ there is a random measure $\,X_{t}^{n}$\thinspace\ in
$\,\mathcal{M}$\thinspace\ satisfying%
\begin{equation}
-\log\mathbf{P}\mathrm{e}^{-\left\langle X_{t}^{n},\varphi\right\rangle
}\;=\;\left\langle \mu,v_{n}(t)_{\!_{\!_{\,}}}\right\rangle ,\qquad\varphi
\in\Phi.
\end{equation}
In fact, $\,X_{t}^{n}$\thinspace\ arises in the following way: Start at time
$\,0$\thinspace\ with $\,\mu.$\thinspace\ During $\,t\in\lbrack0,\frac{1}%
{n}),$\thinspace\ let $\,X_{t}^{n}$\thinspace\ evolve independently in each
point according to continuous state branching with index $\,1+\beta
,$\thinspace\ starting from $\,\mu,$\thinspace\ that is, related to definition
(\ref{def-vn-expl}). Then, at time $\,\frac{1}{n}$\thinspace\ let
$\,X_{1/n-}^{n}$\thinspace\ jump to $\,S_{1/n}^{\alpha}X_{1/n\,-}%
^{n}=:X_{1/n\,}^{n}.$\thinspace\ Then continue with continuous state branching
as before and up to time $\,\frac{2}{n}-,$\thinspace\ etc., until
$\,t$\thinspace\ is reached.

\textrm{F}rom the convergence $\,v_{n}(t)\rightarrow v(t)$\thinspace
\ according to Lemma \ref{P.vn-converge} and domination as in (\ref{dom}),
implying $\,v(t)\downarrow0$\thinspace\ as $\,\varphi\downarrow0,$%
\thinspace\ we get the existence of a random measure $\,X_{t}$\thinspace\ in
$\,\mathcal{M}$\thinspace\ such that $\,X_{t}^{n}\rightarrow X_{t}$%
\thinspace\ in law as $\,n\uparrow\infty.$\thinspace\ Since the map
$\,\mu\mapsto\left\langle \mu,V_{t}\varphi\right\rangle $\thinspace\ is
measurable, via $\,\mu\mapsto X_{t}$\thinspace\ we get a probability kernel
$\,Q_{t}$\thinspace\ in $\,\mathcal{M}$.\thinspace\ \textrm{F}rom the
semigroup property of $\,V\!\varphi$\thinspace\ follows that $\,Q=\left\{
Q_{t}:\,t>0\right\}  $\thinspace\ satisfies Chapman-Kolmogorov. Hence,
$\,Q$\thinspace\ is the transition kernel of a Markov process $\,X$%
\thinspace\ in $\,\mathcal{M},$\thinspace\ which is the desired superprocess.
This finishes the proof.\hfill\rule{0.5em}{0.5em}\bigskip

\noindent\emph{Acknowledgment}\quad We would like to thank Sergio Albeverio,
Zdzis\l aw Brze\'{z}\-niak, Werner Kirsch, and Karl-Theodor Sturm for
discussions on the operators $\,\Delta^{\!(\alpha)}.$


\vfill

\hfill{\scriptsize birth12.tex\quad typeset by \LaTeX}\quad
{\scriptsize printed \today}

\vfill


\newcommand{\noopsort}[1]{}

\end{document}